\documentclass[11pt,a4paper]{article}       
\usepackage[affil-it]{authblk}

\usepackage[utf8]{inputenc}
\usepackage{amsfonts, graphicx, amsmath, amssymb, amsthm}
\usepackage{siunitx}
\usepackage{tikz}
	\usetikzlibrary{calc}
\usepackage{subcaption}
\usepackage{appendix}
\usepackage{url}

\newtheorem{thrm}{Theorem}[section]
\newtheorem{lem}[thrm]{Lemma}

\newtheorem{prop}[thrm]{Proposition}

\newtheorem{assumption}[thrm]{Assumption}

\newtheorem{rem}[thrm]{Remark}

\newcommand{\Pmax}{P_{\max}}
\newcommand{\vmax}{v^{\max}}
\newcommand{\wmax}{w^{\max}}
\newcommand{\wmin}{w^{\min}}
\newcommand{\winit}{w_{\text{in}}}
\newcommand{\wfin}{w_{\text{fin}}}
\newcommand{\wrefp}{w^+}
\newcommand{\wrefm}{w^-}

\usepackage{algorithm}
\usepackage{algorithmicx}
\usepackage{algpseudocode}
\floatname{algorithm}{Algorithm} 

\title{\LARGE \bf A convex reformulation for speed planning of a vehicle under the travel time and energy consumption objectives}

\author{Luca Consolini, Mattia Laurini, and Marco Locatelli
\thanks{All authors are with the University of Parma, Department of Engineering and Architecture, 181/A
Parco Area delle Scienze, 43124 Parma, Italy. E-mails:
{\tt \{luca.consolini, mattia.laurini, marco.locatelli\}@unipr.it}.}}%
\date{}

\begin{document}
\allowdisplaybreaks[4]

\maketitle

\begin{abstract}
In this paper we address the speed planning problem for a vehicle along a predefined path.
A weighted sum of two conflicting objectives, energy consumption and travel time, is minimized.
After deriving a non-convex mathematical model of the problem, we prove that the feasible region 
of this problem is a lattice.
Moreover, we introduce a feasibility-based bound-tightening technique which allows us to derive the minimum and maximum element of the lattice, or establish that the feasible region is empty.
We prove the exactness of a convex relaxation of the non-convex problem, obtained by replacing all constraints with the lower and upper bounds for the variables corresponding to the minimum and maximum elements of the lattice, respectively.
After proving some properties of optimal solutions of the convex relaxation, we exploit them to develop a dynamic programming approach returning  an approximate solution to the convex relaxation, and with time complexity $O(n^2)$, where $n$ is the number of points into which the continuous path is discretized.  
\end{abstract}

\section{Introduction}
\label{sec:intro}
In this paper we address the problem of speed planning for a road vehicle on an assigned path in such a way that two distinct and conflicting objectives, namely travel time and energy consumption,
are minimized.
We consider a simplified model of a half-car.
We denote by:
\begin{itemize}
\item $M$ the mass and $v$ the speed of the vehicle;
\item $\alpha(s)$ the slope angle of the road, which is a function of arc-length position $s$;
\item $F_a(t)$ the aerodynamic drag force, which is a function of time $t$;
\item $T_f$ and $T_r$ the tractive (if positive) or braking (if negative) forces of the front and rear wheels; 
\item $R_f$ and $R_r$ the rolling resistances of the front and rear wheels. 
\end{itemize}
We consider the component of the car's acceleration orthogonal to the road negligible.
Then, the following dynamic equations, obtained from the balance of forces in the direction parallel and orthogonal to the road, can be derived:
\begin{equation}
\begin{aligned}
\label{eqn_dyn_1a}
M \dot v(t) =	&\ T_f(t)+T_r(t)-R_f(t)-R_r(t)-F_a(t) - g M \sin \alpha(s(t)) \\
0 =			&\ W_f(t)+W_r(t)-g M \cos \alpha(s(t)).
\end{aligned}
\end{equation}
From the second equation, the overall load is $W_f(t)+W_r(t)=g M \cos \alpha(s(t))$.
The tractive (or braking) forces must satisfy the following constraints
\begin{equation}
  \label{eqn_max_torques_no_slip}
|T_f| \leq \mu W_f,\qquad |T_r| \leq \mu W_r,
\end{equation}
in order to avoid slipping (we assume that the front and rear axles have the same friction coefficient $\mu$).
From the second equation in~\eqref{eqn_dyn_1a}, we have that
\begin{equation}
\label{eqn_no_slip}
|T_f(t)+T_r(t)| \leq g M \mu \cos \alpha(s(t)).
\end{equation}
We assume that the car has four-wheel drive, and that drive and braking forces are optimally balanced. In this case, drive and braking forces can reach the maximum values compatible with conditions~(\ref{eqn_max_torques_no_slip}), and condition~\eqref{eqn_no_slip} is sufficient to prevent slipping on front and rear wheels. In other words, we assume that, if condition~(\ref{eqn_no_slip}) is satisfied, the vehicle controller is always able to provide the required traction or braking forces.
Setting $F(t)= T_f(t)+ T_r(t)$ as the overall tractive and braking force, and $F_r(t)=R_f(t)+R_r(t)$ as the overall rolling resistance, we rewrite the first equation of~\eqref{eqn_dyn_1a} as (see also Table~4-1 in~\cite{rajamani2011vehicle})
\begin{equation}
\label{eqn_dyn}
M \dot v(t)= -F_a(t) + F(t) - g M \sin \alpha(s(t)) - F_r(t).
\end{equation}
The aerodynamic drag is $F_a(t)=\frac{1}{2} \rho A_f c_d v(t)^2= \Gamma v(t)^2$, where $\rho$ is the air density, $A_f$ is the car cross-sectional area, and $c_d$ is the dimensionless drag coefficient.
We model the overall tires rolling resistance with the simplified model
$F_r(t)= (W_f(t)+W_r(t)) c= g M \cos \alpha(s(t)) c$, where $c$ is the dimensionless rolling resistance coefficient. 
The power consumed (if positive) or recovered (if negative) by the vehicle is given by 
\[
P(t)= \max\{\eta F(t),F(t)\} v(t).
\]
According to the value of $\eta \in [0, 1]$, we obtain different power consumption behaviors: in case of an internal combustion engine (ICE), the vehicle is not recovering energy when braking ($\eta = 0$). Otherwise, in case of hybrid or electric engines ($\eta > 0$), the regenerative braking system of the vehicle enables energy recovery during deceleration, which accounts for a negative value for the consumed power, meaning that power is recovered.

Let $\gamma: [0,s_f] \to \mathbb{R}^2$ be a smooth
function, whose image set $\boldsymbol{\gamma}([0,s_f])$ is the path to be
followed. Here, $\gamma(0)$ and $\gamma(s_f)$ are the initial and final car positions.
We assume that $\gamma$ has arc-length parameterization, that is, is such that
\mbox{$(\forall \lambda \in [0,s_f])$,  $\|\gamma(\lambda)\| =1$}. 
In this way, $s_f$ is the path length.
Let $\lambda: [0,t_f] \to
[0,s_f]$ be a differentiable monotone increasing function, that
represents the vehicle's arc-length position along the curve as a function of time, and
let
$v: [0,s_f]\to [0,+\infty[$ be such that $(\forall t \in
[0,t_f])\ \dot \lambda(t)=v(\lambda(t))$. In this way, $v(s)$
is the derivative of the vehicle arc-length position, which corresponds to the norm of its velocity vector at position $s$. The position of the vehicle as a function of time is given by $x:[0,t_f] \to \mathbb{R}^2, \ x(t)=\gamma(\lambda(t))$.\\
We choose the optimization variable $w(s)=\frac{1}{2} v^2(s)$, that is half the squared speed, as a function of position $s$ (see, for instance,~\cite{verscheure09}). Then, $w'(s) = \frac{dw}{ds} = v(s)\frac{dv}{ds} = v(s)\frac{dv}{v(s)dt} = \frac{dv}{dt}=\dot v(t)$, and~\eqref{eqn_dyn} becomes
\begin{equation}
  \label{eqn_dyn_with_hat}
M w'(s)= -\hat F_a(s) + \hat F(s) - g M \sin \alpha(s) - \hat F_r(s),
\end{equation}
where forces $\hat F_a$, $\hat F$, $\hat F_r$ depend on arc-length $s$. That is, they are functions $[0,s_f] \to \mathbb{R}$ such that, for all $t \in [0,t_f]$,
\[
\hat F_a(\lambda(t))=F_a(t),
\hat F(\lambda(t))=F(t),
\hat F_r(\lambda(t))=F_r(t).
\]
To simplify notation, we write $\hat F_a$, $\hat F$, $\hat F_r$ as $F_a$, $F$, $F_r$ (that is, we redefine these as arc-length functions), and rewrite~(\ref{eqn_dyn_with_hat}) as:
\[
M w'(s)= -F_a(s) + F(s) - g M \sin \alpha(s) - F_r(s).
\]
We can then write our optimization problem as follows:
\begin{equation}
  \label{eqn_prog_cont}
\min_w \int_{s=0}^{s_f} \left(\lambda \max\{\eta F(s),F(s)\} + \frac{1}{\sqrt{2w(s)}}\right)ds
\end{equation}
\vspace{-2em}
\begin{align}
\text{s. t.} &														& \nonumber\\
& w(s) \leq \wmax(s)													& \text{for } s \in [0, s_f], \\
& |F(s)| \leq g M \mu \cos \alpha(s)										& \text{for } s \in [0, s_f], \\
& F(s) \leq \frac{\Pmax}{\sqrt{2w(s)}}										& \text{for } s \in [0, s_f], \\
& M w'(s)= -F_a(s) + F(s) +											& \nonumber\\
&\qquad\qquad\enspace- g M \sin \alpha(s) - F_r(s)										& \text{for } s \in [0, s_f], \\
& w(s) \geq 0														& \text{for } s \in [0, s_f], \\
& w(0) = \winit,\quad w(s_{f})=\wfin.										&
\end{align}
Note that:
\begin{itemize}
\item an upper bound $\wmax(s)$ is imposed for half squared speed $w(s)$, depending on the position along the path, since, e.g., the maximum allowed speed is different along a curve and along a straight road;
\item an upper bound $g M \mu \cos \alpha(s)$ is imposed for the tractive (or braking) force;
\item an upper bound $\Pmax$ is imposed for the power;
\item initial and final speed are fixed to $\winit$ and to $\wfin$, respectively (boundary conditions);
\item the objective function is a weighted sum of the two terms: 
\begin{align*}
& \int_{s=0}^{s_f} \max\{\eta F(s),F(s)\}ds\ \mbox{(energy consumption),} \\
& \int_{s=0}^{s_f} \frac{1}{\sqrt{2w(s)}}ds \ \mbox{(travel time)}. 
\end{align*}
Parameter $\lambda\geq 0$ is the weight assigned to the energy consumption term. It represents the time cost (in seconds) per Joule (its dimension is \unit{\second\per\joule}).
\end{itemize}
We discretize the problem by observing variable and input functions at positions $\{0,h,2h,\ldots,(n-1)h\}$, where $h$ is the discretization step and $n \in \mathbb{N}$ is such that $(n-1)h \approx s_f$ is the length of the path.
We also approximate derivatives with finite differences and the integral in the objective function
with the Riemann sum of the intervals.
Then, the resulting discretized version of the previous problem is:
\[
\min_w \sum_{i=1}^{n-1} \left(\lambda \max\{\eta F_i,F_i\} + \frac{1}{\sqrt{2w_i}}\right)h
\]
\vspace{-2em}
\begin{align*}
\textrm{s. t.} &								& \\
& w_i \leq \wmax_i							& \!\!\!\!\!\!\!\text{for } i \in \{1, \ldots, n\}, \\
& |F_i| \leq g M \mu \cos \alpha_i				& \!\!\!\!\!\!\!\text{for } i \in \{1, \ldots, n-1\},\\
& F_i \leq \frac{\Pmax}{\sqrt{2w_i}}				& \!\!\!\!\!\!\!\text{for } i \in \{1, \ldots, n-1\}, \\
& \frac{M}{h} (w_{i+1}-w_i) = - \Gamma w_i+ F_i +	& \\
& \qquad\qquad\qquad\quad\enspace - g M(\sin \alpha_i + c \cos \alpha_i)		& \!\!\!\!\!\!\!\text{for } i \in \{1, \ldots, n-1\}, \\
& w_i \geq 0								& \!\!\!\!\!\!\!\text{for } i \in \{1, \ldots, n\}, \\
& w_1 = \winit,\quad w_n=\wfin.					&
\end{align*}
After setting $f_i=\frac{F_i}{M}$, we can also write the problem as follows ($\gamma=\frac{\Gamma}{M}$):
\begin{equation}
  \label{eqn_prob_disc}
\min_w \sum_{i=1}^{n-1} \left(\lambda M \max\{\eta f_i,f_i\} + \frac{1}{\sqrt{2w_i}}\right)h
\end{equation}
\vspace{-2em}
\begin{align}
\text{s. t.} &								& \nonumber\\
& w_i \leq \wmax_i							& \!\!\text{for } i \!\in\! \{1, \ldots, n\}, \\
& |f_i| \leq g \mu \cos \alpha_i					& \!\!\text{for } i \!\in\! \{1, \ldots, n\!-\!1\}, \\
& f_i \leq \frac{\Pmax}{M\sqrt{2w_i}}				& \!\!\text{for } i \!\in\! \{1, \ldots, n\!-\!1\}, \\
& \frac{1}{h} (w_{i+1}\!-\!w_i) \!=\! -\gamma w_i + f_i+	& \nonumber\\
& \qquad\qquad\qquad\ - g(\sin \alpha_i \!+\! c \cos \alpha_i)	& \!\!\text{for } i \!\in\! \{1, \ldots, n\!-\!1\}, \\
& w_i \geq 0								& \!\!\text{for } i \!\in\! \{1, \ldots, n\}, \\
& w_1 = \winit,\quad w_n=\wfin.					&
\end{align}
In the following sections we will provide some results to solve efficiently the discretized problem.
Next, the solution of the discretized problem can be employed to derive an approximate solution of the continuous-time problem through interpolation. For instance, one can quadratically interpolate the discrete-time solution of~(\ref{eqn_prob_disc}) with the method we discussed in~\cite{Consolini2017converges}, or use other interpolation techniques. One could discuss if the obtained interpolated solution is a solution of the continuous-time problem~(\ref{eqn_prog_cont}). We expect that, adapting the arguments in~\cite{Consolini2017converges}, it is possible to show that, as the discretization step $h$ goes to $0$, the quadratically interpolated solution of the discretized problem converges to a solution of the continuous-time problem~(\ref{eqn_prog_cont}). However, a precise analysis of this convergence would be lengthy, and is not the focus of the present work, which is centered on the discretized problem~(\ref{eqn_prob_disc}).
\begin{rem}
\label{rem:pareto}
We can also view the speed planning problem as a bi-objective optimization problem, where one objective function, say $g_1(w)$, represents travel time, while the other, say $g_2(w)$, represents fuel consumption.
In a bi-objective optimization problem one aims at detecting Pareto optimal solutions, i.e., feasible solutions $w^\star$ which are {\em not} dominated, i.e., for which there does not exist another feasible solution
$w$ such that $g_1(w^\star)\geq g_1(w)$ and $g_2(w^\star)\geq g_2(w)$, where one of the two inequalities is strict (we refer, e.g. to~\cite{eichfelder2021twenty} for a survey about multiobjective optimization).
The set of all Pareto optimal solutions is called {\em Pareto front}.
A widely employed way to compute Pareto optimal solutions is through the so called {\em linear scalarization method}, where Pareto optimal solutions are computed through
the optimization of the single objective function $\eta g_1(w)+(1-\eta)g_2(w)$ for $\eta\in [0,1]$. Thus, solving the previously introduced discretized problem with different values of the parameter $\lambda$ can also be viewed as an application of the linear scalarization method.
\end{rem}

\subsection{Literature Review} 
When $\lambda=0$, the problem becomes a minimum-time one. On the other hand, if $\lambda$ is large, this problem becomes a minimum-energy one. We review the literature both on minimum-time and minimum-energy problems, since they are strongly related.

\subsubsection{Minimum-time speed planning}
Various works focus on the minimum-time control problem for road vehicles on an assigned path.
Some works represent the vehicle speed as a function of time, others represent the vehicle speed as a function of the arc-length position along the assigned path. In some cases, the second choice considerably simplifies the resulting optimization problem.
Among the works that make the first choice,~\cite{MunOllPraSim:94,munoz1998speed} propose an iterative algorithm based on the concatenation of third degree polynomials. Reference~\cite{SolNun2006} achieves a solution with an algorithm based on the five-splines scheme of~\cite{bianco2006optimal}, while~\cite{Villagra-et-al2012} proposes a closed-form speed profile obtained as a concatenation of three different classes of speed profiles. Finally,~\cite{CheHeBuHanZha2014} proposes an algorithm which returns a piecewise linear speed profile.
Among the works that parameterize the speed as a function of the arc-length,~\cite{Minari16,minSCL17} consider a problem that, after finite-element discretization, can be reformulated as a convex one. This allows deriving efficient solution methods such as those presented in~\cite{CLMNV19,minSCL17,LippBoyd2014}. In particular, reference~\cite{minSCL17} presents an algorithm, with linear-time computational complexity with respect to the number of variables, that provides an optimal solution after spatial discretization. Namely, the path is divided into $n$ intervals of equal length, and the problem is approximated with a finite dimensional one in which the derivative of $v$ is substituted with a finite difference approximation. Some other works directly solve this problem in continuous-time~\cite{consolini2020solution,Frego16,Velenis2008}. In particular, paper~\cite{consolini2020solution} computes directly the exact continuous-time solution without performing a finite-dimensional reduction, proposing a method which is very simple and can be implemented very efficiently, so that it is particularly well suited for real-time speed planning applications. Reference~\cite{RaiPerCGL17jerk} solves the problem with an additional jerk constraint through a heuristic approach that computes a speed profile by bisection; the proposed method is very efficient, however, the optimality of the obtained profile is not guaranteed.

\subsubsection{Energy consumption minimization}
Other works seek the speed profile that minimizes energy consumption on an assigned path.
Various works optimize the energy utilization of metro vehicles for a given trip.
For instance,~\cite{wang2011optimal} formulates the optimal speed profile for energy saving problem as a Mixed Integer Linear Programming (MILP) one.
In~\cite{kang2011ga}, Kang et al. present an algorithm for optimizing a train speed profile by controlling the coasting point, using Genetic Algorithms (GA).
Reference~\cite{calderaro2014algorithm} finds the sequence of basic control regimes that minimizes the energy consumed on a given path, taking into account the track topology (slopes and curves), the mechanical characteristics of the vehicle, the electrical characteristics of the feeding line, and the effect of regenerative braking. 

 Other works minimize fuel consumption for private vehicles.
However, human drivers are not aware of the optimal velocity profile for a given
 route. Indeed, the globally optimal velocity trajectory depends on
 many factors, and its calculation requires intensive computations~\cite{gustafsson2009automotive,russell2002integrated}. With advancements in communication, sensors, and in-vehicle computing, real-time optimizations are becoming more feasible. For example, some public transportation vehicles in Europe now communicate with traffic lights~\cite{koenders2008cooperative}, and the USA is experimenting with broadcasting red light timings for safety~\cite{intersections2008cooperative}. Several algorithms have been proposed for optimizing speed profiles, such as predicting optimal velocity profiles for approaching traffic lights~\cite{asadi2010predictive}, or minimizing energy use for a given route with a traffic light~\cite{ozatay2012analytical}. Despite these developments, many algorithms still require expensive on-board systems and have limited real-time applicability. To address these limitations, cloud computing offers a promising solution for real-time speed optimization~\cite{wollaeger2012cloud}. Reference~\cite{ozatay_cloud-based_2014} builds on that by implementing a real-time Speed Advisory System (SAS), that uses cloud computing to generate optimal velocity profiles based on traffic and geographic data, offering a global optimization approach. 
 A more aggressive approach is rapidly accelerating to a
 given speed, followed by a period of coasting to a lower speed,
 which is called the ``pulse-and-glide'' (PnG) strategy~\cite{li2015effect,li2016fuel,li2015mechanism,xu2015fuel}. PnG driving can
 achieve significant fuel savings in vehicles with continuously
 variable transmissions (CVT)~\cite{li2015mechanism} and step-gear mechanical
 transmissions~\cite{xu2015fuel}. The optimality and the increasing efficiency are proven by theoretical calculations. Reference~\cite{kim2019real} introduces a real-time implementable pulse-and-glide algorithm in the speed-acceleration domain, achieving 3\%--5\% fuel savings without compromising optimality. Some papers focus on the study of road slope~\cite{hellstrom2009look,kamal2011ecological,musardo2005ecms,xu2018design}, which significantly impacts fuel economy, requiring optimization of speed and control systems. Model predictive control (MPC) methods, like those by Kamal et al.~\cite{kamal2011ecological} and Hellstr\"{o}m et al.~\cite{hellstrom2009look}, improve fuel efficiency, but face high computational demands due to system non-linearity. Non-predictive methods, such as ECMS~\cite{musardo2005ecms}, offer reduced computational loads and near-optimal fuel savings by relying on instantaneous slope data. Both approaches are viable for connected automated vehicles (CAVs), but require balancing fuel economy and computational efficiency for practical implementation. Some references (for instance, ~\cite{zheng2020study}) do consider a multi-objective problem for route selection, but not for speed planning.
We did not find any work that formulates speed planning as a multi-objective minimization problem, where the objectives are traveling time and energy consumption.
Anyway, some works do consider strictly related problems. For instance, some references optimize only energy consumption, but introduce a constraint on the traveling time. For instance,~\cite{calderaro2014algorithm} uses this approach for optimizing the speed profile of a metro vehicle, whose longitudinal dynamics are similar to a road vehicle. Paper~\cite{xu2018design} minimizes only energy consumption, but bounds the minimum velocity, which limits traveling time from below. However, these references use an approach quite different from ours. They do not find a convex relaxation of the original problem, and do not find the globally optimal solution. Namely, reference~\cite{calderaro2014algorithm} proposes a two-step optimization, while~\cite{xu2018design} uses a generic non-linear solver.

\subsection{Statement of contribution}
The present work builds upon and extends our previous work on speed planning~\cite{ARDIZZONI2025}.
While the earlier formulation already addressed the trade-off between travel time and energy consumption, it did not incorporate a constraint on final speed. 
With respect to that work, in this paper we provide the following novel contributions:
\begin{itemize}
\item we explicitly include a constraint on final speed, which is required in some applications (think about the case when the vehicle must stop at the end of path, i.e., its final speed must be null);
\item we develop a feasibility-based bound-tightening procedure that enables us to derive a convex reformulation of the original non-convex problem, different from the one discussed in~\cite{ARDIZZONI2025} and, actually, simpler since expressed solely through box constraints;
\item the feasibility-based bound-tightening procedure also allows us establishing a necessary and sufficient condition for the feasibility of the problem, providing a clear criterion to verify feasibility;
\item exploiting structural properties of the resulting convex formulation, we introduce a Dynamic Programming (DP) algorithm that provides an approximate solution of the original problem with a computational complexity of $O(n^2)$, where $n$ is the number of discretization steps. This method is tailored to handle large scale instances efficiently, while preserving high solution quality;
\item finally, through a set of computational experiments, we show that our approach achieves computation times approximately three orders of magnitude faster than those obtained with a commercial solver for convex problems, while maintaining a very small relative error on the objective value.
\end{itemize}

\subsection{Outline of the paper}
After having observed that the optimization problem to be solved is a non-convex one, in Section~\ref{sec:feasboundtighten} we introduce a feasibility-based bound-tightening technique.
In Section~\ref{sec:lattice} we prove that the feasible region of our problem is a lattice and from that we derive a necessary and sufficient condition to establish its non-emptiness, also discussing a procedure to verify such condition, based on the iterated application of the bound-tightening technique discussed in Section~\ref{sec:feasboundtighten}. In Section~\ref{sec:convrel} we prove that a convex relaxation of the non-convex problem, based on the outcome of the bound-tightening technique, is exact (i.e., its solution is equivalent to the optimal solution of the non-convex problem).
In Section~\ref{sec:optsolprop}, we prove some properties of optimal solutions of the convex relaxation (and, thus, of the non-convex problem).
In Section~\ref{sec:dynprog} we exploit the properties proved in Section~\ref{sec:optsolprop} to develop a DP approach, with time complexity $O(n^2)$, which returns an approximate solution of the convex relaxation.
We show that the distance between the speeds explored by the DP approach and those of the optimal solution is bounded from above by a quantity proportional to the discretization step $h$.
In Section~\ref{sec:compexp} we present some computational experiments.
Finally, in Section~\ref{sec:concl} we draw the conclusions and discuss a possible extension to more general force and power constraints.

\section{Feasibility-based bound-tightening technique}
\label{sec:feasboundtighten}
The constraints of our problem are the following: for $j \in \{1, \ldots, n-1\}$, $i \in \{1, \ldots, n\}$,
\begin{align}
& \frac{1}{h} (w_{j+1}-w_{j}) + \gamma w_{j} + g (\sin \alpha_{j} + c \cos \alpha_{j})\leq \frac{\Pmax}{M\sqrt{2w_{j}}}		& \label{eq:powermax} \\
& \frac{1}{h} (w_{j+1}-w_{j}) + \gamma w_{j} + g (\sin \alpha_{j} + c \cos \alpha_{j})\leq g\mu\cos \alpha_{j}		& \label{eq:forcemax1} \\
& \frac{1}{h} (w_{j+1}-w_{j}) + \gamma w_{j} + g (\sin \alpha_{j} + c \cos \alpha_{j})\geq -g\mu\cos \alpha_{j}	& \label{eq:forcemax2} \\
& w_i \leq \wmax_i					& \label{eq:speedmax} \\
& w_i \geq 0						& \label{eq:nonneg} \\
& w_1 = \winit, \qquad w_n=\wfin.		& \label{eq:limitcond}
\end{align}
We rewrite constraints~\eqref{eq:speedmax}--\eqref{eq:limitcond} as follows:
\begin{align}
& w_i \leq \wmax_i & i \in \{1, \ldots, n\} \label{eq:wmax} \\
& w_i \geq \wmin_i & i \in \{1, \ldots, n\} \label{eq:wmin},\!\!
\end{align}
where $\wmin_i=0$, for $i\in \{2,\ldots,n-1\}$, while $\wmin_1=\wmax_1=\winit$ and $\wmin_n=\wmax_n=\wfin$ (note that the two boundary conditions~\eqref{eq:limitcond} are split into constraints $\winit\leq w_1\leq \winit$ and $\wfin\leq w_n\leq \wfin$).
We denote by $X$ the set of points $w$, $w\in \mathbb{R}^n$, fulfilling constraints~\eqref{eq:powermax}--\eqref{eq:limitcond}.
Such region is non-convex due to the non-convexity of power constraints~\eqref{eq:powermax}.
Then, our problem is the following (non-convex) one
\begin{equation}
\label{eq:probfix}
\begin{array}{lll}
\min\limits_{w\in X, f} & \sum\limits_{i=1}^{n-1} h \left(\lambda M\max\{\eta f_i,f_i\} \!+\! \frac{1}{\sqrt{2w_i}}\right) \\ 
& f_i=\frac{1}{h} (w_{i+1}-w_i) + \gamma w_i + g (\sin \alpha_i + c\cos \alpha_i),	\enspace i \in \{1,\ldots,n\!-\!1\}.
\end{array}
\end{equation}
While convex problems are usually solvable in polynomial time (see, e.g.,~\cite{nesterov1994interior}), non-convex problems are hard to solve. Nevertheless, some non-convex problems have the so called
{\em hidden convexity property}, i.e., they can be reformulated into equivalent convex problems. The best known optimization problem with such property is the trust region one, where a non-convex quadratic function is minimized over the unit sphere (see~\cite{rendl1997semidefinite}). We refer to~\cite{xia2020survey} for a survey on hidden convex optimization. In this paper we will show that problem~\eqref{eq:probfix} has the hidden convexity property, and, thus, can be efficiently solved.
Now, let us introduce the following functions, $i \in \{1,\ldots,n-1\}$:
\begin{equation}
\label{eq:defell}
\ell_i(w)=(1-h\gamma)w+h\min\left\{\frac{\Pmax}{M\sqrt{2w}},g\mu \cos \alpha_i\right\}.
\end{equation}
After defining the critical half squared speeds $\hat{w}_i$, $i \in \{1,\ldots,n-1\}$, as the solution of
\begin{equation}\label{eq:critical_speed}
\frac{\Pmax}{M\sqrt{2\hat{w}_i}} = g\mu \cos \alpha_i,
\end{equation}
we have that:
\[
\ell_i(w)=
\begin{cases}
(1-h\gamma)w+hg\mu\cos \alpha_i & w\leq \hat{w}_i \\
(1-h\gamma)w+h\frac{\Pmax}{M\sqrt{2w}} & w> \hat{w}_i .
\end{cases}
\]
Next, we introduce the following assumption.
\begin{assumption}
\label{assum:3}
It holds that for each $i\in\{1,\ldots,n-1\}$:
\begin{equation}
\label{ineq:incrass}
1 - h\gamma -h\frac{\Pmax}{M\left(2\hat{w}_i\right)^{\frac{3}{2}}} \geq 0.
\end{equation}
\end{assumption}
Note that this assumption is fulfilled if the discretization step $h$ is small enough.
Under this assumption we can prove the following lemma.
\begin{lem}
\label{lem:incr}
Under Assumption~\ref{assum:3}, functions $\ell_i$, $i=1,\ldots,n-1$, are increasing for $w>0$.
\end{lem}
\begin{proof}
Under Assumption~\ref{assum:3} it holds that $1-h\gamma>0$ so that $\ell_i$ is increasing for $w\leq \hat{w}_i$. For $w>\hat{w}_i$ we have that the derivative of $\ell_i$ is:
\[
\ell'_i(w)=1-h\gamma -h\frac{\Pmax}{M\left(2w\right)^{\frac{3}{2}}}.
\]
Since $\ell'_i$ is increasing with respect to $w$, if $\ell'_i(\hat{w}_i)\geq 0$, then $\ell'_i(w)>0$ for all $w>\hat{w}_i$ (i.e., $\ell_i$ is increasing for $w>\hat{w}_i$).
Now it is enough to observe that
\[
\ell'_i(\hat{w}_i)=1 - h\gamma -h\frac{\Pmax}{M\left(2\hat{w}_i\right)^{\frac{3}{2}}},
\]
so that, under Assumption~\ref{assum:3}, $\ell_i$ is increasing also for $w>\hat{w}_i$.
\end{proof}
Now, let us consider region $X_1\supset X$, defined as follows:
\begin{equation}
\label{eq:powforcons}
\begin{aligned}
X_1=&\left\{ w \mid \frac{1}{h} (w_{i+1}-w_i) + \gamma w_i + g (\sin \alpha_i + c\cos \alpha_i)\leq\right. \\
&\leq\left.\min\left\{\frac{\Pmax}{M\sqrt{2w_{i}}},g\mu\cos \alpha_i \right\},\ i \in \{1,\ldots,n-1\},\ w\geq \wmin \right\},
\end{aligned}
\end{equation}
that is, the region defined by maximum power and maximum force constraints~\eqref{eq:powermax} and~\eqref{eq:forcemax1}, and by lower limit constraints~\eqref{eq:wmin}.

\noindent
For a fixed value $\bar{w}_{i+1}$, $i\in \{1,\ldots,n-1\}$, let us consider the $i$-th maximum power and maximum force constraints and the $i$-th lower limit constraint:
\begin{equation}
\label{eq:powforcons1}
\begin{aligned}
 \frac{1}{h} (\bar{w}_{i+1}-w_i) + \gamma w_i + g (\sin \alpha_i + c\cos \alpha_i)	& \leq\frac{\Pmax}{M\sqrt{2w_i}} \\
 \frac{1}{h} (\bar{w}_{i+1}-w_i) + \gamma w_i + g (\sin \alpha_i + c\cos \alpha_i)	& \leq g \mu\cos \alpha_i \\
 \wmin_i												& \leq w_i.
\end{aligned}
\end{equation}
 This system of inequalities can be rewritten as:
\begin{equation}
\label{eq:powforcons2}
\begin{aligned}
\bar{w}_{i+1} + hg (\sin \alpha_i + c\cos \alpha_i)	& \leq \ell_i(w_i) \\
 \wmin_i						& \leq w_i.
\end{aligned}
\end{equation}
We denote by $\xi_1^{\wmin}(\bar{w}_{i+1})$ the smallest value which can be assigned to $w_i$ so that
the $i$-th maximum power, maximum force and lower limit constraints are fulfilled (i.e,~\eqref{eq:powforcons2} is fulfilled).
We prove the following lemma.
\begin{lem}
\label{lem:xi}
We have that $\xi_1^{\wmin}(\bar{w}_{i+1})=\wmin_i$ if $\ell_i(\wmin_i)\geq \bar{w}_{i+1} + hg (\sin \alpha_i + c\cos \alpha_i)$, while for $\ell_i(\wmin_i)< \bar{w}_{i+1} + hg (\sin \alpha_i + c\cos \alpha_i)$, we have that
$\xi_1^{\wmin}(\bar{w}_{i+1})>\wmin_i$ and $\xi_1^{\wmin}$ is increasing with $\bar{w}_{i+1}$. 
\end{lem} 
\begin{proof}
In view of Lemma~\ref{lem:incr}, $\ell_i$ is increasing for $w>0$. Then, if $\ell_i(\wmin_i)\geq \bar{w}_{i+1} + hg (\sin \alpha_i + c\cos \alpha_i)$, we have that $\xi_1^{\wmin}(\bar{w}_{i+1})=\wmin_i$, while
for $\ell_i(\wmin_i)< \bar{w}_{i+1} + hg (\sin \alpha_i + c\cos \alpha_i)$, there exists a unique value $\xi_1^{\wmin}(\bar{w}_{i+1})>\wmin_i$ which can be assigned to $w_i$ such that
$\ell_i(\xi_1^{\wmin}(\bar{w}_{i+1}))= \bar{w}_{i+1} + hg (\sin \alpha_i + c\cos \alpha_i)$.
Moreover, the value
$\xi_1^{\wmin}(\bar{w}_{i+1})$ is increasing with $\bar{w}_{i+1}$. 
\end{proof}
Now we define a point $l$ as follows:
\begin{equation}
\label{eq:recurtildew}
l_n=\wfin,\ \ l_{j}=\xi_1^{\wmin}(l_{j+1})\qquad j\in \{1,\ldots,n-1\}.
\end{equation}
The following proposition shows that each value $l_{j}$ is the smallest value which can be attained by variable $w_{j}$ at points $w\in X_1$.
\begin{prop}
\label{prop:smallest}
For $j \in \{1,\ldots,n\}$, if $x \in X_1$, then $w_{j} \geq l_{j}$.
\end{prop}
\begin{proof}
We prove this by induction. The result is obviously true for $j=n$ since in $X_1$ we have that $w_n$ is bounded from below by $\wfin$.
Now, for $j< n$, let us assume that the result is true for all $i\in \{j+1,\ldots,n\}$ and let us prove it for $j$. In particular, the inductive assumption guarantees that
$l_{j+1}$ is the smallest value of variable $w_{j+1}$ which can be assigned to such variable at points in $X_1$.
Moreover, for each $w \in X_1$ for which $w_{j+1}=\bar{y}$, where $\bar{y}$ is a fixed value, the smallest value which can be attained by $w_{j}$ is $\xi_1^{\wmin}(\bar{y})$, that is, the smallest solution
of the system of inequalities~\eqref{eq:powforcons1} or the equivalent system~\eqref{eq:powforcons2} for $\bar{w}_{j+1}=\bar{y}$.
In view of Lemma~\ref{lem:xi}, the minimum value for $\xi_1^{\wmin}(\bar{y})$ is attained when $\bar{y}$ is smallest, that is, when $\bar{y}=l_{j+1}$, since, by the inductive assumption, $l_{j+1}$ is the smallest possible value of $w_{j+1}$ at points in $X_1$. Thus, the smallest possible value for $w_{j}$ at points in $X_1$ is $\xi_1^{\wmin}(l_{j+1})=l_{j}$, as we wanted to prove. 
\end{proof}
Since $X_1 \supset X$, we can also conclude that values $l_j$, $j \in \{1,\ldots,n\}$, are valid lower bounds for the variables $w_j$ at all points in $X$ (i.e., all feasible solutions of Problem~\eqref{eq:probfix}).
In other words, we can say that recursive relation~\eqref{eq:recurtildew} defines a {\em feasibility-based bound-tightening} procedure, that is, a procedure that employs some of the constraints defining region $X$, namely maximum power and maximum force constraints~\eqref{eq:powermax} and~\eqref{eq:forcemax1}, together with lower limit constraints~\eqref{eq:wmin}, to tighten the (lower) bounds on the variables (for an introduction to bound-tightening techniques we refer, e.g., to~\cite{belotti2009branching,schichl2005interval,tawarmalani2004global}).

Now, we can introduce further feasibility-based bound-tightening procedures.
Let us introduce set $X_2\supset X$ defined as follows: 
\[
\begin{aligned}
X_2=&\left\{w \mid \frac{1}{h} (w_{i+1}-w_i) + \gamma w_i + g (\sin \alpha_i + c\cos \alpha_i)\leq\right.\\
&\leq\left.\min\left\{\frac{\Pmax}{M\sqrt{2w_{i}}},g\mu \cos \alpha_i\right\},\ i \in \{1,\ldots,n-1\},\ 
w\leq \wmax\right\},
\end{aligned}
\]
that is, the region defined by maximum power and maximum force constraints~\eqref{eq:powermax} and~\eqref{eq:forcemax1}, and by upper limit constraints~\eqref{eq:wmax}.

For a fixed value $\bar{w}_{i}$, $i\in \{1,\ldots,n-1\}$, let us consider the $i$-th maximum power and maximum force constraints and the $(i+1)$-th maximum speed constraint:
\begin{align*}
 \frac{1}{h} (w_{i+1}-\bar{w}_{i}) + \gamma \bar{w}_{i} + g (\sin \alpha_i + c\cos \alpha_i)	& \leq\frac{\Pmax}{M\sqrt{2\bar{w}_i}} \\
 \frac{1}{h} (w_{i+1}-\bar{w}_{i}) + \gamma \bar{w}_{i} + g (\sin \alpha_i + c\cos \alpha_i)	& \leq g \mu \cos \alpha_i\\
 w_{i+1}														& \leq \wmax_{i+1}.
\end{align*}
 Then, a valid upper bound for variable $w_{i+1}$ is:
\[
w_{i+1}\leq \xi_2^{\wmax}(\bar{w}_i)=\min\left\{\wmax_{i+1},\ell_i(\bar{w}_i)-h g (\sin \alpha_i + c\cos \alpha_i)\right\}.
\]
Similarly to the proof of Proposition~\ref{prop:smallest}, we can show that the recursive relation
\[
u_1=\winit,\ \ u_{j+1}=\xi_2^{\wmax}(u_j)\qquad j\in \{1,\ldots,n-1\},
\]
defines a set of valid upper bounds for variables $w_j$.

\noindent
Next, let us introduce set $X_3\supset X$ defined as follows: 
\begin{align*}
X_3&=\left\{w \mid -\frac{1}{h} (w_{i+1}\!-\!w_i) \!-\! \gamma w_i \!-\! g (\sin \alpha_i \!+\! c\cos \alpha_i)\leq \right.\\
&\left. \phantom{\frac11}\leq g\mu\cos \alpha_i,\ \ i \in \{1,\ldots,n-1\},\ w\geq \wmin \right\},
\end{align*}
that is, the region defined by minimum force constraints~\eqref{eq:forcemax2},
and by lower limit constraints~\eqref{eq:wmin}.

\noindent
For a fixed value $\bar{w}_{i}$, $i\in \{1,\ldots,n-1\}$, let us consider the $i$-th minimum force constraint and the $(i+1)$-th lower limit constraint:
\begin{align*}
 -\frac{1}{h} (w_{i+1}-\bar{w}_{i}) - \gamma \bar{w}_{i} - g (\sin \alpha_i + c\cos \alpha_i)	& \leq g \mu \cos \alpha_i\\
 \wmin_{i+1}													& \leq w_{i+1} .
\end{align*}
 Then, a valid lower bound for variable $w_{i+1}$ is:
\begin{align*}
 w_{i+1}\geq \xi_3^{\wmin}(\bar{w}_i)=&\max\left\{\wmin_{i+1},(1\!-\!h\gamma)\bar{w}_i \!-\!h g (\sin \alpha_i \!+\! c\cos \alpha_i+\mu\cos \alpha_i)\right\}.
\end{align*}
Similarly to the proof of Proposition~\ref{prop:smallest}, we can show that the recursive relation
\[
l'_1=\winit,\ \ l'_{j+1}=\xi_3^{\wmin}(l_j)\qquad j\in \{1,\ldots,n-1\},
\]
defines a set of valid lower bounds for variables $w_j$. 

\noindent
Finally, let us introduce set $X_4\supset X$ defined as follows: 
\begin{align*}
X_4&=\left\{ w \mid -\frac{1}{h} (w_{i+1}\!-\!w_i) \!-\! \gamma w_i \!-\! g (\sin \alpha_i \!+\! c\cos \alpha_i)\leq\right. \\
&\left.\phantom{\frac11}\leq g\mu\cos \alpha_i,\ i \in \{1,\ldots,n-1\},\ w\leq \wmax\right\},
\end{align*}
that is, the region defined by minimum force constraints~\eqref{eq:forcemax2}, 
and by upper limit constraints~\eqref{eq:wmax}.
For a fixed value $\bar{w}_{i+1}$, $i\in \{1,\ldots,n-1\}$, let us consider the $i$-th minimum force constraint and the $i$-th maximum speed constraint:
\begin{align*}
 -\frac{1}{h} (\bar{w}_{i+1}-w_{i}) - \gamma w_{i} - g (\sin \alpha_i + c\cos \alpha_i)	& \leq g \mu\cos \alpha_i \\
 w_i														& \leq \wmax_{i} .
\end{align*}
Then, a valid upper bound for variable $w_{i}$ is:
\begin{align*}
w_{i}\leq& \xi_4^{\wmax}(\bar{w}_{i+1})=\min\left\{\wmax_i, \frac{\bar{w}_{i+1}+h g (\sin \alpha_i + c\cos \alpha_i+\mu\cos \alpha_i)}{1-h\gamma}\right\}.
\end{align*}
Similarly to the proof of Proposition~\ref{prop:smallest}, we can show that the recursive relation
\[
u'_n=\wfin,\ \ u'_{j}=\xi_4^{\wmax}(u'_{j+1})\qquad j\in \{1,\ldots,n-1\},
\]
defines a set of valid upper bounds for the variables $w_j$. 

\noindent
Note that $X=X_1\cap X_2\cap X_3\cap X_4$.

\section{Lattice structure and a necessary and sufficient condition for feasibility}
\label{sec:lattice}
We denote by $\vee$ and $\wedge$ the component-wise maximum and minimum of two vectors of dimension $n$, respectively, that is, given $w,w'\in \mathbb{R}^n$, we have that for each $i \in \{1,\ldots,n\}$:
\[
[w\vee w']_i=\max\{w_i,w'_i\},\qquad [w\wedge w']_i=\min\{w_i,w'_i\}.
\]
\begin{prop}
\label{prop:lattice}
Let Assumption~\ref{assum:3} hold.
If $w, w'\in X$, then $w \vee w'\in X$ and $w \wedge w'\in X$, that is, $(X,\vee,\wedge)$ is a lattice. 
\end{prop}
\begin{proof}
We notice that $w, w'\in X$ implies that $w\vee w'\in X$ and $w\wedge w'\in X$ fulfill lower limit constraints~\eqref{eq:wmin}, upper limit constraints~\eqref{eq:wmax}, and boundary conditions~\eqref{eq:limitcond}. Therefore, we only need to prove that $w\vee w'\in X$ and $w\wedge w'\in X$ also fulfill maximum power constraints~\eqref{eq:powermax}, and maximum and minimum force constraints~\eqref{eq:forcemax1}--\eqref{eq:forcemax2}.
Also recalling definition~\eqref{eq:defell} of $\ell_i$, we can rewrite these constraints as follows:
for $i \in \{1,\ldots,n-1\}$
\begin{equation}
\label{eq:reform}
\begin{aligned}
&w_{i+1} \leq \ell_i(w_i)-h g (\sin \alpha_i + c\cos \alpha_i) \\
&w_{i+1} \geq (1-h\gamma)w_i -h g [\sin \alpha_i + (c+\mu)\cos \alpha_i],
\end{aligned}
\end{equation}
where the first inequality is equivalent to the $i$-th maximum power constraint~\eqref{eq:powermax} and maximum force constraint~\eqref{eq:forcemax1}.
Let us assume, w.l.o.g., that $w_{i+1}\wedge w'_{i+1}=w_{i+1}$ and $w_{i+1}\vee w'_{i+1}=w'_{i+1}$. 

\noindent
If $w_i\wedge w'_i=w_i$ and $w_i\vee w'_i=w'_i$, then constraints~\eqref{eq:reform} are obviously fulfilled by 
$w\wedge w'$ and by $w\vee w'$. 

\noindent
Therefore, let us see what happens when $w_i\wedge w'_i=w'_i$ and $w_i\vee w'_i=w_i$. 
First of all, we observe that:
\begin{align*}
w_{i+1} \!\wedge\! w'_{i+1} &= w_{i+1}\leq w'_{i+1} \leq \ell_i(\underbrace{w'_i}_{w_{i}\wedge w'_{i}})-h g (\sin\alpha_i + c\cos \alpha_i) \\
w_{i+1} \!\wedge\! w'_{i+1} & =w_{i+1}\!\geq\! (1\!-\!h\gamma)w_i \!-\!h g [\sin\alpha_i \!+\! (c\!+\!\mu)\!\cos \alpha_i] \\
&\geq (1\!-\!h\gamma) \!\!\underbrace{w'_{i}}_{w_{i}\wedge w'_{i}}\!\! -\!h g [\sin\alpha_i \!+\! (c\!+\!\mu)\cos \alpha_i],
\end{align*}
that is, constraints~\eqref{eq:reform} are fulfilled by $w\wedge w'$. 
Next, recalling the increasingness of $\ell_i$ established in Lemma~\ref{lem:incr} under Assumption~\ref{assum:3}, we observe that:
\begin{align*}
w_{i+1}\vee w'_{i+1} & =w'_{i+1}\leq \ell_i(w'_i)-h g (\sin\alpha_i + c\cos \alpha_i)\leq\\
&\leq \ell_i(\underbrace{w_i}_{w_{i}\vee w'_{i}})-h g (\sin\alpha_i + c\cos \alpha_i) \\
w_{i+1}\vee w'_{i+1} & =w'_{i+1}\geq w_{i+1}\geq (1-h\gamma)\!\!\underbrace{w_i}_{w_{i}\vee w'_{i}}\!\! -h g [\sin\alpha_i \!+\! (c\!+\!\mu)\cos \alpha_i],
\end{align*}
that is, constraints~\eqref{eq:reform} are fulfilled by $w\vee w'$. 
\end{proof}
\begin{rem}
There are other problems whose feasibility region is a lattice. That is, given two feasible solutions, their component-wise maximum and minimum are still feasible. For instance, we used this property in our previous works~\cite{consolini2020solution,minSCL17}. Also, lattice structures have been leveraged in different optimization contexts. For instance,~\cite{Topkis1978} shows that if a function is submodular, its minima are a sublattice.
\end{rem}

Since $(X,\vee,\wedge)$ is a lattice, $X\neq \varnothing$ implies that the component-wise maximum and minimum over $X$, that is, points $z',y'$ defined as follows:
\[
z'_i=\max_{w\in X} w_i,\qquad y'_i=\min_{w\in X} w_i, \qquad i \in \{1,\ldots,n\},
\]
are such that $z',y'\in X$. But how can we establish whether $X\neq \varnothing$ and, in case it is, how can we compute $z',y'$?
We proceed as follows. Let
\[
z_i=\max_{w\in X_2\cap X_4} w_i,\quad y_i=\min_{w\in X_1\cap X_3} w_i, \quad i \in \{1,\ldots,n\},
\]
that is:
\begin{itemize}
\item $z_i$ is the maximum value which can be attained by $w_i$ if we ignore the lower bound constraints on the variables, that is, the non-negativity constraints for variables $w_i$, $i\in \{2,\ldots,n-1\}$ together with $w_1\geq \winit$
and $w_n\geq \wfin$;
\item $y_i$ is the minimum value which can be attained by $w_i$ if we ignore the maximum speed constraints (where $\wmax_1=\winit$ and $\wmax_n=\wfin$).
\end{itemize}
Note that $X_1\cap X_3$ and $X_2\cap X_4$ are both non-empty lattices, so that $z\in X_1\cap X_3$ and $y\in X_2\cap X_4$ (later on we will introduce a procedure for their computation).
We prove the following proposition stating a necessary and sufficient condition for $X\neq \varnothing$.
\begin{prop}
\label{prop:empty}
We have that $X\neq \varnothing$ if and only if $y\leq z$. Moreover, if $X\neq \varnothing$, then $z'=z$ and $y'=y$.
\end{prop}
\begin{proof}
If $y\leq z$, then $\wmax_i\geq z_i\geq y_i\geq 0$, for $i\in \{2,\ldots,n-1\}$, and:
\[
\begin{array}{ll}
\wmax_1=\winit\geq z_1\geq y_1\geq \winit & \implies z_1=y_1=\winit \\
\wmax_n=\wfin\geq z_n\geq y_n\geq \wfin & \implies z_n=y_n=\wfin,
\end{array}
\]
so that $y,z\in X$ and $X\neq \varnothing$. 
\noindent
Moreover, since $X\subset X_1\cap X_3, X_2\cap X_4$, in this case, $z$ and $y$ are the component-wise maximum and minimum over $X$, respectively (i.e., $z'=z$ and $y'=y$).

\noindent
Instead, if $y\not\leq z$, we can conclude that $X=\varnothing$. Indeed, since each $z_i$ is an upper limit for $w_i$ over $X_2\cap X_4$ and, thus, over $X$, while each $y_i$ is a lower limit for $w_i$ over $X_1\cap X_3$ and, thus, over $X$,
if for some $i$ it holds that $y_i>z_i$, we can conclude that $X=\varnothing$.
\end{proof}
The next question is how to compute $z,y$. We proceed as follows. First, we introduce the following functions:
for $j \in \{1,\ldots,n-1\}$,
\begin{align}
&B_1(u)=p	\qquad\mbox{where } p_1=u_1,\ p_{j+1}=\xi_2^l(p_{j}),	& \label{eq:B_1}\\
&B_2(u)=p	\qquad\mbox{where } p_n=u_n,\ p_j=\xi_4^u(p_{j+1}),	& \label{eq:B_2}\\
&B_3(l)\; =p	\qquad\mbox{where } p_n=l_n,\ p_j=\xi_1^l(p_{j+1}),		& \label{eq:B_3}\\
&B_4(l)\; =p	\qquad\mbox{where } p_1=l_1,\ p_{j+1}=\xi_3^l(p_{j}),	& \label{eq:B_4}
\end{align}
where we recall that functions $\xi_1,\xi_2,\xi_3,\xi_4$ have been defined in Section~\ref{sec:feasboundtighten}.
Then, we can apply Algorithm~\ref{alg:buildsol} to compute $z,y$. More precisely, we show that Algorithm~\ref{alg:buildsol} is an iterative algorithm converging to $z,y$ or establishing that $X=\varnothing$.
The algorithm works as follows. At line~\ref{lin:1} vectors $u^1$ and $l^1$ are initialized with vectors $\wmax$ and $\wmin$, respectively.
Next, at iteration $k$ we compute: 
\begin{itemize}
\item at line~\ref{lin:2} a new upper bound vector $u^{k+\frac{1}{2}}$ for variables $w$. By~\eqref{eq:B_1}, it holds that $u^{k+\frac{1}{2}}\leq u^k$;
\item at line~\ref{lin:3} a new upper bound vector $u^{k+1}$ for variables $w$. By~\eqref{eq:B_2}, it holds that $u^{k+1}\leq u^{k+\frac{1}{2}}$;
\item at line~\ref{lin:4} a new lower bound vector $l^{k+\frac{1}{2}}$ for variables $w$. By~\eqref{eq:B_3}, it holds that $l^{k+\frac{1}{2}}\geq l^k$;
\item at line~\ref{lin:5} a new lower bound vector $l^{k+1}$ for variables $w$. By~\eqref{eq:B_4}, it holds that $l^{k+1}\geq l^{k+\frac{1}{2}}$.
\end{itemize}
Then, we check whether the distance between $u_k$ and $u_{k+1}$, and the one between $l^k$ and $l^{k+1}$ is within a given small threshold distance $\varepsilon>0$. If so, we stop
and return the vectors $u^{k+1}$ and $l^{k+1}$. Otherwise, if $u^{k+1}\not\geq l^{k+1}$, since $z\leq u^{k+1}$ and $l^{k+1}\leq y$, we have that $z\not\geq y$ and, in view of Proposition~\ref{prop:empty}, we can stop
establishing that $X=\varnothing$.

\noindent
The following proposition proves a convergence result for Algorithm~\ref{alg:buildsol}.
\begin{prop}
If we set $\varepsilon=0$, then Algorithm~\ref{alg:buildsol} either stops after a finite number of iterations returning $z,y$, or generates two infinite sequences $\{u^k\}$ and $\{l^k\}$, the former converging from above to $z$, the latter converging from below to $y$.
\end{prop}
\begin{proof}
First we observe that the two sequences $\{u^k\}$ and $\{l^k\}$ are non-increasing and non-decreasing, respectively. Moreover, for each $k$ it holds that $u^k\geq u^{k+\frac{1}{2}}\geq u^{k+1}\geq z$ and $l^k\leq l^{k+\frac{1}{2}}\leq l^{k+1}\leq y$.
By definition of $u^{k+\frac{1}{2}}, u^{k+1}, l^{k+\frac{1}{2}}, l^{k+1}$, we have that, for $i\in \{1,\ldots,n-1\}$:
\begin{equation}
\label{eq:convcons}
\begin{array}{l}
\frac{1}{h} (u^{k+\frac{1}{2}}_{i+1}-u^{k+\frac{1}{2}}_{i}) + \gamma u^{k+\frac{1}{2}}_{i} + g (\sin \alpha_i + c\cos \alpha_i) \leq \frac{\Pmax}{M\sqrt{2u^{k+\frac{1}{2}}_{i}}} \\
\frac{1}{h} (u^{k+\frac{1}{2}}_{i+1}-u^{k+\frac{1}{2}}_{i}) + \gamma u^{k+\frac{1}{2}}_{i} + g (\sin \alpha_i + c\cos \alpha_i) \leq g \mu\cos \alpha_i \\
\frac{1}{h} (u^{k+1}_{i+1}-u^{k+1}_{i}) + \gamma u^{k+1}_{i}+ g (\sin \alpha_i + c\cos \alpha_i)					 \geq  -g \mu \cos \alpha_i \\
\frac{1}{h} (l^{k+\frac{1}{2}}_{i+1}-l^{k+\frac{1}{2}}_{i}) + \gamma l^{k+\frac{1}{2}}_{i} + g (\sin \alpha_i + c\cos \alpha_i)	 \leq \frac{\Pmax}{M\sqrt{2l^{k+\frac{1}{2}}_{i}}} \\
\frac{1}{h} (l^{k+\frac{1}{2}}_{i+1}-l^{k+\frac{1}{2}}_{i}) + \gamma l^{k+\frac{1}{2}}_{i} + g (\sin \alpha_i + c\cos \alpha_i)	 \leq g \mu \cos \alpha_i \\
\frac{1}{h} (l^{k+1}_{i+1}-l^{k+1}_{i}) + \gamma l^{k+1}_{i}+ g (\sin \alpha_i + c\cos \alpha_i)						 \geq -g \mu\cos \alpha_i.
\end{array}
\end{equation}
If at iteration $k$ we have that $\|u^{k+1}-u^k\|=\|l^{k+1}-l^k\|=0$, then $u^k=u^{k+\frac{1}{2}}=u^{k+1}$ and $l^k=l^{k+\frac{1}{2}}=l^{k+1}$, and, in view of~\eqref{eq:convcons}, $u^k\in X_2\cap X_4$, while $l^k\in X_1\cap X_3$, so that $z=u^k$ and $y=l^k$.

\noindent
Instead, if $\|u^{k+1}-u^k\|=\|l^{k+1}-l^k\|=0$ never occurs, then we observe that sequences $\{u^k\}$ and $\{l^k\}$ are monotonic and, if we never stop, they are also bounded, since $u_k$ cannot fall below $\wmin$, while $l^k$ cannot be larger than $\wmax$. Therefore, the sequences converge to $\bar{u}\geq z$ and $\bar{l}\leq y$, respectively. Taking the limit for $k\rightarrow \infty$ in~\eqref{eq:convcons}, we have that $\bar{u}\in X_2\cap X_4$, while $\bar{l}\in X_1\cap X_3$, so that $z=\bar{u}$ and $y=\bar{l}$. 
\end{proof}

\begin{algorithm}
\caption{\label{alg:buildsol} Procedure to compute the vectors $z,y$ or to establish that $X=\varnothing$}
\begin{algorithmic}[1]
\Procedure{ComputeZY}{$\wmax,\wmin,\varepsilon>0$}
\State Set $u^1=\wmax$, $l^1=\wmin$, $k=1$, $stop={\tt false}$ \label{lin:1}\;
\While{{\tt true}}
\State Set $u^{k+\frac{1}{2}}=B_1(u^k)$ \label{lin:2}\;
\State Set $u^{k+1}=B_2(u^{k+\frac{1}{2}})$ \label{lin:3}\;
\State Set $l^{k+\frac{1}{2}}=B_3(l^k)$\label{lin:4}\;
\State Set $l^{k+1}=B_4(l^{k+\frac{1}{2}})$\label{lin:5}\;
\If{$\|u^{k+1}-u^k\|\leq \varepsilon\ \ \mbox{{\tt and}}\ \ \|l^{k+1}-l^k\|\leq \varepsilon$ \label{lin:6}}
\State \Return{$u^{k+1}, l^{k+1}$}\;
\ElsIf{$u^{k+1}\not\geq l^{k+1}$\label{lin:7}}
\State{$X=\varnothing$}
\Else
\State Set $k=k+1$\;
\EndIf
\EndWhile
\EndProcedure
\end{algorithmic}
\end{algorithm}

\section{An exact convex relaxation}
\label{sec:convrel}
In this section we employ the component-wise minimum $y\in X$ and component-wise maximum $z\in X$ of the feasible region $X$ of Problem~\eqref{eq:probfix}, computed through
Algorithm~\ref{alg:buildsol} , to define a convex relaxation of that problem, and we will prove that such convex relaxation is exact, that is, it shares the same optimal solution and optimal value of Problem~\eqref{eq:probfix}.

First of all we introduce two reference speeds:
\begin{equation*}
\wrefp=\left(2\lambda M \gamma\right)^{-\frac{2}{3}},\qquad \wrefm=\left(2\eta \lambda M \gamma\right)^{-\frac{2}{3}}.
\end{equation*}
Note that these speeds are such that they minimize the $i$-th component of the objective function of~\eqref{eq:probfix}, assuming that the speed is kept constant at step $i$ and $i+1$. More specifically, if we express the $i$-th term of the objective function as a function of $w$, and recall that $\frac{1}{h}(w_{i+1} - w_i) = 0$ (since $w_{i+1} = w_i$), we obtain the following two cases
\[\
\begin{cases}
h\!\left(\lambda M \gamma w_i \!+\! g (\sin \alpha_i \!+\! c\cos \alpha_i) \!+\! \frac{1}{\sqrt{2w_i}}\right)\!, &\!\!\!\!\mbox{if } \gamma w_i \!+\! g (\sin \alpha_i + c\cos \alpha_i) \geq 0.\\
h\!\left(\lambda M \eta(\gamma w_i \!+\! g (\sin \alpha_i \!+\! c\cos \alpha_i)) \!+\! \frac{1}{\sqrt{2w_i}}\right)\!, &\!\!\!\!\mbox{otherwise.}
\end{cases}
\]
Observe that, by computing the first derivative of the two previous terms, we have that the minimum of the first and the second case are attained at $\wrefp$ and $\wrefm$, respectively.

Next, we notice that the objective function of Problem~\eqref{eq:probfix}, without constant factor $h$, can be rewritten as follows:
\begin{gather*}
\sum_{i=1}^{n-1} \left[\lambda M \max\{\eta f_i,f_i\} + \frac{1}{\sqrt{2w_i}}\right]=\\
=\sum_{i=1}^{n-1} \left[\lambda M \eta f_i +(1-\eta)\lambda M \max\{0,f_i\} + \frac{1}{\sqrt{2w_i}}\right].	
\end{gather*}
By definition of $f_i$, we have that
\[
\sum_{i=1}^{n-1} f_i = \frac{\wfin-\winit}{h}+g\!\sum_{i=1}^{n-1}(\sin(\alpha_i) + c\cos(\alpha_i))+\gamma \!\sum_{i=1}^{n-1} w_i,
\]
Then, if we introduce the function:
\begin{align*}
F(w,f)=&\sum_{j=2}^{n-1}\left[\eta \lambda M \gamma w_j +\frac{1}{\sqrt{2w_j}}\right]
+(1-\eta)\lambda M \sum_{j=1}^{n-1} \max\{0,f_j\},
\end{align*}
 we have that the objective function of Problem~\eqref{eq:probfix} can be rewritten as:
\[
 h F(w,f) + h\left[\frac{\wfin-\winit}{h}+ g\sum_{i=1}^{n-1}(\sin(\alpha_i)+c\cos(\alpha_i))\right].
\]
In what follows, we will also rewrite function $F(w,f)$ as the sum $F_1(w)+F_2(f)$, where:
\[
F_1(w)=\sum_{j=2}^{n-1}\left[\eta \lambda M \gamma w_j +\frac{1}{\sqrt{2w_j}}\right],\qquad
F_2(f)=(1-\eta)\lambda M \sum_{j=1}^{n-1} \max\{0,f_j\}.
\]
As previously mentioned, let $y,z$ be the lower and upper bound vectors for variables $w$ computed via the bound-tightening procedure Algorithm~\ref{alg:buildsol}, respectively. Consider the following convex problem, obtained by removing all force and power constraints, and introducing the lower and upper limits $y_i$ and $z_i$ for variables $w_i$, $i \in \{1,\ldots,n\}$:
\begin{equation}
\label{eq:newconvrel}
\begin{aligned}
\min_{w,f}\	& h \left[F(w,f) + \frac{\wfin-\winit}{h}+g\sum_{i=1}^{n-1}(\sin(\alpha_i) + c\cos(\alpha_i))\right]\\
& f_i=\frac{1}{h} (w_{i+1}-w_i) +\gamma w_i+ g(\sin \alpha_i + c\cos(\alpha_i))\ \text{for } i \in \{1, \ldots, n-1\} \\
& y_i\leq w_i \leq z_i		\enspace\quad\qquad\qquad\qquad\qquad\qquad\qquad\qquad\, \text{for } i \in \{1, \ldots, n\}.
\end{aligned}
\end{equation}
This problem is a convex relaxation of Problem~\eqref{eq:probfix}. Note that, by disregarding the constant addend and factor, we can replace the objective function of 
Problem~\eqref{eq:newconvrel} with $F(w,f)$.
We prove that, under suitable assumptions, an optimal solution of Problem~\eqref{eq:newconvrel} is also feasible and, thus, optimal for~\eqref{eq:probfix}, thus proving that convex relaxation~\eqref{eq:newconvrel} is actually a convex reformulation of Problem~\eqref{eq:probfix}.
\noindent
Now, for a given feasible solution $w$ of Problem~\eqref{eq:newconvrel}, we denote by $(w^{j,\pm}(\delta), f^{j,\pm}(\delta))$, for some $\delta>0$, the solution defined as follows:
\[
\begin{array}{rll}
(\forall r\neq j)		& w_r^{j,\pm}(\delta) =w_r\\ 
				& w_j^{j,\pm}(\delta) =w_j\pm\delta\\ 
(\forall r\neq j-1,\ j)	& f_{r}^{j,\pm}(\delta) =f_{r}\\ 
				& f_{j-1}^{j,\pm}(\delta) =f_{j-1}\pm\frac{\delta}{h} \\ 
				& f_{j}^{j,\pm}(\delta) =f_{j}\mp\frac{\delta}{h}\pm \gamma \delta.
\end{array}
\]
We prove two lemmas, returning conditions for the speed at the beginning and at the end of a sequence of consecutive steps where an optimal solution of~\eqref{eq:newconvrel} violates maximum force/power constraints
(Lemma~\ref{lem:maxforce}) and minimum force constraints (Lemma~\ref{lem:minforce}).
\begin{lem}
\label{lem:maxforce}
Let $w$ be an optimal solution of Problem~\eqref{eq:newconvrel}.
Consider a maximal sequence $i,\ldots,k$, $k\geq i$, of consecutive steps where the force violates the maximum force or power value, that is, 
\[
\begin{array}{lllll}
f_j		&>\min\left\{g\mu\cos(\alpha_j),\frac{\Pmax}{M\sqrt{2w_j}}\right\} & j \in \{i,\ldots,k\} \vspace{1pt}\\
f_{k+1}	&\leq \min\left\{g\mu\cos(\alpha_{k+1}),\frac{\Pmax}{M\sqrt{2w_{k+1}}}\right\}& \mbox{ \tt or } \enspace k+1=n \vspace{1pt}\\ 
f_{i-1}	&\leq \min\left\{g\mu\cos(\alpha_{i-1}),\frac{\Pmax}{M\sqrt{2w_{i-1}}}\right\}& \mbox{ \tt or } \enspace i=1.
\end{array}
\]
Then:
\[
\begin{array}{lll}
w_i= z_i\enspace 				&\mbox{ \tt or } & \enspace w_i\geq \wrefp \\
w_{k+1}\leq \wrefp\enspace	&\mbox{ \tt or } & \enspace w_{k+1}=y_{k+1}.
\end{array}
\]
\end{lem}
\begin{proof}
Let us first prove that either $w_{k+1}\leq \wrefp$ or $w_{k+1}=y_{k+1}$. If $w_{k+1}=y_{k+1}$ we are done. Otherwise, 
we first notice that solution $(w^{k+1,-}(\delta), f^{k+1,-}(\delta))$ is feasible for Problem~\eqref{eq:newconvrel}, and that 
\[
F(w^{k+1,-}(\delta), f^{k+1,-}(\delta))-F(w,f)\leq -\lambda M \gamma \delta +\frac{\delta}{2(w_{k+1})^\frac{3}{2}},
\]
 where the left-hand side is negative for $w_{k+1}>\wrefp$, thus contradicting the optimality of $w$.

\noindent
Similarly, if $w_i=z_i$, we are done. Instead, if $w_i<z_i$, then solution $(w^{i,+}(\delta), f^{i,+}(\delta))$ is feasible for Problem~\eqref{eq:newconvrel}, and 
\[
F(w^{i,+}(\delta), f^{i,+}(\delta))-F(w,f)\leq \lambda M \gamma \delta -\frac{\delta}{2(w_{i})^\frac{3}{2}},
\]
 where the left-hand side is negative for $w_{i}<\wrefp$, thus contradicting the optimality of $w$.
\end{proof}
Analogously to the proof of Lemma~\ref{lem:maxforce}, we can prove the following lemma.
\begin{lem}
\label{lem:minforce}
Consider a maximal sequence $i,\ldots,k$, $k\geq i$, of consecutive steps where the minimum force constraint is violated, that is, 
\[
\begin{array}{llll}
f_j		& <-g\mu\cos(\alpha_j)							&	& j \in \{i,\ldots,k\} \\  
f_{k+1}	& \geq-g\mu\cos(\alpha_{k+1})&\enspace \mbox{ \tt or }\enspace k+1=n	& \\  
f_{i-1}	& \geq -g\mu\cos(\alpha_{i-1})&\enspace \mbox{ \tt or }\enspace i=1.	&
\end{array}
\]
Then:
\[
\begin{array}{lll}
w_{k+1}=z_{k+1}\enspace		& \mbox{ \tt or }& \enspace w_{k+1}\geq \wrefm \\  
w_i\leq \wrefm\enspace	& \mbox{ \tt or }& \enspace w_i=y_i.
\end{array}
\]
\end{lem}
Now, let us introduce the following assumptions.
\begin{assumption}
\label{ass:condaccdec}
It holds that:
\begin{equation}
1-h\gamma(1 + \lambda P_{\max})	 > 0 \label{cond1}
\end{equation}
and $(\forall i \in \{1, \ldots, n\})$
\begin{align}
z_i\geq \wrefp \implies& \min\left\{g\mu\cos(\alpha_i),\frac{P_{\max}}{M\sqrt{2\wrefp}}\right\} + \label{cond2}\\
&-\gamma \wrefp-g(\sin(\alpha_i) + c\cos(\alpha_i)) > 0  \nonumber
\end{align}
\begin{equation}
-g[\sin(\alpha_i) + (c+\mu)\cos(\alpha_i)] - \gamma \wrefm < 0. \label{cond3}
\end{equation}
\end{assumption}
The following remark discusses why the assumptions above are mild ones.
\begin{rem}
First, we notice that all the conditions of Assumption~\ref{ass:condaccdec} can be checked in advance. Moreover, such conditions are mild ones. More in detail:
\begin{itemize}
\item Condition~\eqref{cond1} is always fulfilled provided that we take the discretization step $h$ small enough. 
\item
Condition~\eqref{cond2} is not related to the discretization step but depends on physical properties of the vehicle and of the road. It is equivalent to require that we are still able to accelerate at the first reference speed $\wrefp$.
Observing that the conditions need to be checked only when $z_i\geq \wrefp$, and that the left-hand side is decreasing with $\wrefp$, the smallest possible value for the left-hand side is
\[
 \min\left\{g\mu\cos(\alpha_i),\frac{P_{\max}}{M\sqrt{2z_i}}\right\} -\gamma z_i-g(\sin(\alpha_i) + c\cos(\alpha_i)). 
\]
We considered three kinds of vehicles: a Fiat 500e full electric (whose specifications are reported in Table~\ref{tab:vehicle_spec}), an ICE Fiat Panda, and a Citro\"{e}n Ami full electric. In order to test the robustness of Condition~\eqref{cond2} under conservative assumptions, for all three vehicles we replaced the curb mass with their Gross Vehicle Weight (GVW). This choice reduces the power-to-mass ratio ($51.2$, $35.4$, and $8$~\unit{\watt\per\kilogram}, respectively), and therefore makes it harder for the vehicle to accelerate on steep grades. In other words, using GVW instead of curb mass acts as a ``stress test'' for the feasibility of Condition~\eqref{cond2}. We also set $\mu = 0.4$, corresponding to wet asphalt and thus modeling non-ideal adhesion conditions. Moreover, for the Fiat Panda and the Citro\"en Ami we considered $\gamma$ equal to $2.8681\cdot10^{-4}$ and $6.3733\cdot10^{-4}$, respectively, and $c$ equal to $0.011$ and $0.013$, respectively. For each vehicle, Figure~\ref{fig:slopes} shows the relation between maximum speed (the values along the $x$-axis correspond to $\sqrt{2z_i}$
measured in  {\em km/h}) and maximum slope for which Condition~\eqref{cond2} is satisfied (the values along the $y$-axis). Note that for the first two vehicles the maximum speed is capped at $\SI{130}{\kilo\meter\per\hour}$, which is the speed limit in Italy, while for the Citro\"{e}n Ami the maximum attainable speed is $\SI{45}{\kilo\meter\per\hour}$.
\begin{figure}[h!]
	\centering
	\includegraphics[width=.6\columnwidth]{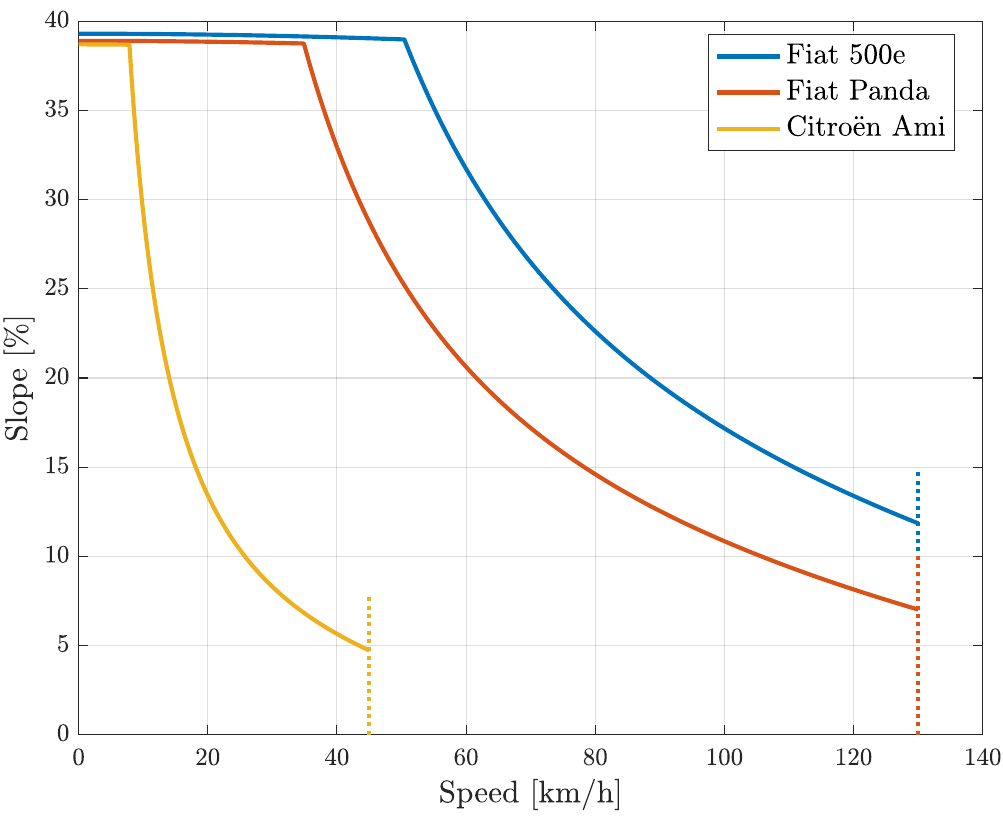}
	\caption{Maximum slopes that satisfy Condition~\eqref{cond2} achievable by each vehicle for different speed values.}
	\label{fig:slopes}
\end{figure}
We notice that as the maximum speed decreases, the condition is fulfilled up to extremely high slopes, but also at the largest speed of the vehicles ($\SI{130}{\kilo\meter\per\hour}$ for the first two vehicles and $\SI{45}{\kilo\meter\per\hour}$ for the third one) violations only occur at relatively high slopes.
\item Condition~\eqref{cond3} also has a physical interpretation and is equivalent to require that we are still able to decelerate at the second reference speed $\wrefm$.
It is more easily fulfilled with respect to Condition~\eqref{cond2}. Indeed, we notice that, independently from $\wrefm$, it is certainly fulfilled if $\tan(\alpha_i)\geq -(\mu+c)$. 
Again referring to a vehicle with the specifications reported in Table~\ref{tab:vehicle_spec} (Fiat 500e full electric) and to a very slippery road ($\mu=0.1$), the condition is fulfilled for $\alpha_i\geq -6^\circ$ (equivalent to a $10.5\%$ slope), while over dry asphalt (for which we can take $\mu=0.7$), the condition is fulfilled for $\alpha_i \geq -35^\circ$ (equivalent to a $70\%$ slope!)
\end{itemize}
Therefore, we can observe that for common road vehicles and non-extreme slopes, the conditions introduced in Assumption~\ref{ass:condaccdec} are mild ones.
\end{rem}
Now, we prove the following theorem.
\begin{thrm}
\label{thm:main}
Let Assumption~\ref{ass:condaccdec} hold.
Then, an optimal solution $w$ of Problem~\eqref{eq:newconvrel} is feasible and, thus, optimal for Problem~\eqref{eq:probfix}.
\end{thrm}
\begin{proof}
Let us assume, by contradiction, that there exists a maximal sequence $i,\ldots,k$ of steps where the maximum force and/or the power constraint are violated. Then, according to Lemma~\ref{lem:maxforce}, one of the following
holds:
\begin{itemize}
\item $w_i=z_i$;
\item $w_{k+1}=y_{k+1}$;
\item $w_i\geq \wrefp,\ w_{k+1}\le \wrefp$.
\end{itemize}
Let us first consider the case $w_i=z_i$. Since
\[
f_i> \min\left\{g\mu\cos(\alpha_i),\frac{P_{\max}}{M\sqrt{2w_i}}\right\},
\]
also recalling definition~\eqref{eq:defell} of function $\ell$, we must have that
\begin{align*}
w_{i+1}\!>& (1-h\gamma) z_i - hg(\sin(\alpha_i) + c\cos(\alpha_i)) + h\min\left\{g\mu\cos(\alpha_i),\frac{P_{\max}}{M\sqrt{2z_i}}\right\} \!=\\
=&\ell_i(z_i)- hg(\sin(\alpha_i) + c\cos(\alpha_i))\geq z_{i+1},
\end{align*}
where the second inequality follows from the fact that $z$ is feasible for~\eqref{eq:probfix}.
But $w_{i+1}>z_{i+1}$ is not possible since it contradicts the feasibility of $w$ for Problem~\eqref{eq:newconvrel}.

\noindent
Next, let us consider the case $w_{k+1}=y_{k+1}$. Since
\[
f_k> \min\left\{g\mu\cos(\alpha_k),\frac{P_{\max}}{M\sqrt{2w_k}}\right\},
\]
recalling again definition~\eqref{eq:defell} of function $\ell_k$, we have that
\[
\ell_k(w_{k})< y_{k+1} + hg(\sin(\alpha_k) + c\cos(\alpha_k))\leq \ell_k(y_k),
\] 
where the second inequality follows from the fact that $y$ is feasible for~\eqref{eq:probfix}. Then,
due to the fact that $\ell_k$ is increasing, we have that $w_k<y_k$, which is not possible since it contradicts the feasibility of $w$ for Problem~\eqref{eq:newconvrel}. 
\noindent
Finally, let us assume that $w_i\geq \wrefp,\ w_{k+1}\le \wrefp$. We only discuss the case in which $\wrefp\geq \hat{w}_i$. The other case is proved in a completely analogous way.
Let $w_i$ be equal to $\wrefp+\Delta$, for some $\Delta\geq 0$. We would like to prove that $w_{i+1}> \wrefp$. We have that $f_i> \frac{P_{\max}}{M\sqrt{2w_i}}$ implies
\[
w_{i+1}> \wrefp+\Delta+h\left(\frac{P_{\max}}{M\sqrt{2\wrefp+\Delta}}-\gamma \wrefp-\gamma\Delta -g\sin(\alpha_i) -gc\cos(\alpha_i)\right).
\]
Therefore, to prove that $w_{i+1}> \wrefp$, it is enough to prove that
\[
\wrefp+\Delta+h\left(\frac{P_{\max}}{M\sqrt{2\wrefp+\Delta}}-\gamma \wrefp-\gamma\Delta -g\sin(\alpha_i) -gc\cos(\alpha_i)\right) > \wrefp,
\]
or, equivalently, that for all $i$ such that $z_i\geq \wrefp$:
\begin{equation}
\label{eq:nonsotto1}
\begin{aligned}
(1-h\gamma)\Delta &+h\frac{P_{\max}}{M\sqrt{2\wrefp+\Delta}}-h\gamma \wrefp -hg(\sin(\alpha_i) + c\cos(\alpha_i)) > 0.
\end{aligned}
\end{equation}
Once we have proved that $w_{i+1}>\wrefp$, by induction, we can also prove that $w_{i+2},\ldots,w_k,w_{k+1}>\wrefp$, thus leading to a contradiction. 
The left-hand side of~\eqref{eq:nonsotto1} is increasing at some $\Delta\geq 0$ if: 
\[
1-h\gamma -h\frac{P_{\max}}{2M(\wrefp+\Delta)^\frac{3}{2}}.
\]
Therefore, if
\[
1-h\gamma -h\frac{P_{\max}}{2M(\wrefp)^\frac{3}{2}}=1-h\gamma(1 + \lambda P_{\max}) >0,
\]
which is true in view of~\eqref{cond1}, then
the left-hand side of~\eqref{eq:nonsotto1} is increasing 
for all $\Delta\geq 0$ and attains its minimum at $\Delta=0$. 
Then, if for all $i$ such that $z_i\geq \wrefp$ it holds that
\[
h\left[\frac{P_{\max}}{M\sqrt{2\wrefp}} -\gamma \wrefp-g(\sin(\alpha_i) + c\cos(\alpha_i))\right]>0,
\]
which is guaranteed by~\eqref{cond2}, 
we can conclude that the left-hand side of~\eqref{eq:nonsotto1} is always strictly positive, which concludes the proof.

\noindent
In a completely similar way we can prove that minimum force constraints are never violated.
\end{proof}

\section{Properties of optimal solutions}
\label{sec:optsolprop}
In this section we prove some properties of optimal solutions of Problem~\eqref{eq:probfix} or, equivalently, of its exact convex relaxation~\eqref{eq:newconvrel}.
First of all, we prove a proposition which gives indication about speed $w_j$ under different possible combinations of the signs of forces $f_{j-1},f_j$.
\begin{prop}
\label{prop:combposs}
Assume that
\begin{equation}
\label{eq:asscomb}
-\frac{(1-\eta)}{h}+ \gamma < 0.
\end{equation}
Let $(w,f)$ be an optimal solution of Problem~\eqref{eq:newconvrel}. Then, the following results hold:
\[
\begin{array}{|l|l|l|l|}
\hline
\mbox{Case}	& \multicolumn{2}{|c|}{\mbox{Combination}}	& \mbox{Speed} \\ \hline \hline
{\bf a)}		&f_{j-1}<0	& f_j<0	& w_j=y_j\geq \wrefm\ \mbox{\tt or}\ w_j=z_j\leq \wrefm \mbox{\tt or}\ w_j=\wrefm \\ \hline
{\bf b)}		&f_{j-1}<0	& f_j>0	& w_j=z_j \\ \hline
{\bf c)}		&f_{j-1}>0	& f_j<0	& w_j=y_j \\ \hline
{\bf d)}		&f_{j-1}>0	& f_j>0	& w_j=y_j\geq \wrefp\ \mbox{\tt or}\ w_j=z_j\leq \wrefp \mbox{\tt or}\ w_j=\wrefp \\ \hline
{\bf e)}		&f_{j-1}<0	& f_j=0	& w_j\geq \wrefm\ \mbox{\tt or}\ w_j=z_j \\ \hline
{\bf f)}		&f_{j-1}>0	& f_j=0	& w_j=y_j\ \mbox{\tt or}\ w_j\leq \wrefp \\ \hline
{\bf g)}		&f_{j-1}=0	& f_j<0	& w_j=y_j\ \mbox{\tt or}\ w_j\leq \wrefm \\ \hline
{\bf h)}		&f_{j-1}=0	& f_j>0	& w_j\geq \wrefp \ \mbox{\tt or}\ w_j=z_j \\ \hline
\end{array}
\]
\end{prop}
\begin{proof}
We first prove the cases {\bf a)}-{\bf b)} and {\bf d)}-{\bf h)}, while case {\bf c)} will be the last to be proved.

{\bf a)} We notice that in case $w_j\in (y_j,z_j)$, $(w^{j,\pm}(\delta), f^{j,\pm}(\delta))$ are both feasible,
\[
F(w^{j,\pm}(\delta), f^{j,\pm}(\delta))-F(w,f)\approx \pm \eta \lambda M \gamma \delta \mp \frac{\delta}{2(w_{j})^\frac{3}{2}},
\]
and at least one of the two values on the right-hand side is lower than 0, thus contradicting the optimality of $(w,f)$, unless $w_j=\wrefm$.
If $w_j=z_j$, then $(w^{j,-}(\delta), f^{j,-}(\delta))$ is feasible, and 
\[
F(w^{j,-}(\delta), f^{j,-}(\delta))-F(w,f)\approx -\eta \lambda M \gamma \delta + \frac{\delta}{2(w_{j})^\frac{3}{2}},
\]
where the right-hand side is negative if $w_j>\wrefm$, thus contradicting the optimality of $(w,f)$. Finally, if
$w_j=y_j$, then $(w^{j,+}(\delta), f^{j,+}(\delta))$ is feasible, and 
\[
F(w^{j,+}(\delta), f^{j,+}(\delta))-F(w,f)\approx \eta \lambda M \gamma \delta - \frac{\delta}{2(w_{j})^\frac{3}{2}},
\]
where the right-hand side is negative if $w_j<\wrefm$, thus contradicting the optimality of $(w,f)$. 

\noindent
{\bf b)} We notice that in case $w_j<z_j$, $(w^{j,+}(\delta), f^{j,+}(\delta))$ is feasible and that
\begin{align*}
F(w^{j,+}(\delta),& f^{j,+}(\delta))-F(w,f) \approx -\lambda(1-\eta)M\frac{\delta}{h}+\lambda M \gamma \delta - \frac{\delta}{2(w_{j})^\frac{3}{2}}<0,
\end{align*}
where the last inequality follows from~\eqref{eq:asscomb}, so that optimality of $(w,f)$ is contradicted.

\noindent
{\bf d)} We notice that in case $w_j\in (y_j,z_j)$, $(w^{j,\pm}(\delta), f^{j,\pm}(\delta))$ are both feasible, that
\[
F(w^{j,\pm}(\delta), f^{j,\pm}(\delta))-F(w,f)\approx \pm \lambda M \gamma \delta \mp \frac{\delta}{2(w_{j})^\frac{3}{2}},
\]
and at least one of the two values on the right-hand side is lower than 0, thus contradicting the optimality of $(w,f)$, unless $w_j=\wrefp$.
In case $w_j=z_j$, then $(w^{j,-}(\delta), f^{j,-}(\delta))$ is feasible and
\[
F(w^{j,-}(\delta), f^{j,-}(\delta))-F(w,f)\approx -\lambda M \gamma \delta + \frac{\delta}{2(w_{j})^\frac{3}{2}},
\]
where the right-hand side is negative if $w_j>\wrefp$, thus contradicting the optimality of $(w,f)$. Finally, if
$w_j=y_j$, then $(w^{j,+}(\delta), f^{j,+}(\delta))$ is feasible and
\[
F(w^{j,+}(\delta), f^{j,+}(\delta))-F(w,f)\approx \lambda M \gamma \delta - \frac{\delta}{2(w_{j})^\frac{3}{2}},
\]
where the right-hand side is negative if $w_j<\wrefp$, thus contradicting the optimality of $(w,f)$. 

\noindent
{\bf e)} We notice that in case $w_j<z_j$, $(w^{j,+}(\delta), f^{j,+}(\delta))$ is feasible, and that
\[
F(w^{j,+}(\delta), f^{j,+}(\delta))-F(w,f)\approx \eta \lambda M \gamma \delta - \frac{\delta}{2(w_{j})^\frac{3}{2}},
\]
where the right-hand side is negative if $w_j<\wrefm$, thus contradicting the optimality of $(w,f)$.

\noindent
{\bf f)} We notice that, if $w_j>y_j$, then $(w^{j,-}(\delta), f^{j,-}(\delta))$ is feasible, and 
\[
F(w^{j,-}(\delta), f^{j,-}(\delta))-F(w,f)\approx -\lambda M \gamma \delta + \frac{\delta}{2(w_{j})^\frac{3}{2}},
\]
where the right-hand side is negative if $w_j>\wrefp$, thus contradicting the optimality of $(w,f)$.

\noindent
{\bf g)} We notice that, if $w_j>y_j$, then $(w^{j,-}(\delta), f^{j,-}(\delta))$ is feasible, and 
\[
F(w^{j,-}(\delta), f^{j,-}(\delta))-F(w,f)\approx -\eta \lambda M \gamma \delta + \frac{\delta}{2(w_{j})^\frac{3}{2}},
\]
where the right-hand side is negative if $w_j>\wrefm$, thus contradicting the optimality of $(w,f)$.

\noindent
{\bf h)} We notice that in case $w_j<z_j$, $(w^{j,+}(\delta), f^{j,+}(\delta))$ is feasible, and that
\[
F(w^{j,+}(\delta), f^{j,+}(\delta))-F(w,f)\approx \lambda M \gamma \delta - \frac{\delta}{2(w_{j})^\frac{3}{2}},
\]
where the right-hand side is negative if $w_j<\wrefp$, thus contradicting the optimality of $(w,f)$.

The proof of case {\bf c)} requires a bit more work.

{\bf c)} If $w_j=y_j$, then we are done. Therefore, let us assume that $w_j>y_j$. For mathematical reasons,
we need to introduce a third reference speed
\[
w^{{\tt ref}, \eta}= \left[\frac{1}{2\left(\eta \lambda M \gamma+\frac{\lambda M (1-\eta)}{h}\right)}\right]^\frac{2}{3},
\]
which, for $\eta<1$, converges to 0 as $h\rightarrow 0$. 
Since $w_j>y_j$, we have that $(w^{j,-}(\delta), f^{j,-}(\delta))$ is feasible and 
\begin{align*}
F(w^{j,-}(\delta)&, f^{j,-}(\delta))-F(w,f)\approx -\lambda(1-\eta)M\frac{\delta}{h}-\lambda \eta M \gamma \delta + \frac{\delta}{2(w_{j})^\frac{3}{2}},
\end{align*}
where the right-hand side is negative when
$w_j> w^{{\tt ref}, \eta}$, so that the optimality of $(w,f)$ is contradicted. Then, $w_j\leq w^{{\tt ref}, \eta}$ must hold.
Moreover, if 
$w^{{\tt ref}, \eta}<z_j$,
which certainly holds for $h$ small enough because $w^{{\tt ref}, \eta}$ converges to 0 as $h\rightarrow 0$,
we have that also $(w^{j,+}(\delta), f^{j,+}(\delta))$ is feasible and that
\begin{align*}
F(w^{j,+}(\delta),& f^{j,+}(\delta))-F(w,f)\approx \lambda(1-\eta)M\frac{\delta}{h}+\lambda \eta M \gamma \delta - \frac{\delta}{2(w_{j})^\frac{3}{2}},
\end{align*}
where the right-hand side is negative when
$w_j< w^{{\tt ref}, \eta}$, so that the optimality of $(w,f)$ is contradicted. Then, $w_j\geq w^{{\tt ref}, \eta}$ must hold, which, combined with $w_j\leq w^{{\tt ref}, \eta}$, leads to $w_j=w^{{\tt ref}, \eta}$ as the only possible alternative to $w_j=y_j$. Therefore, to complete the proof of case {\bf c)}, we need to prove that $w_j=w^{{\tt ref}, \eta}$ cannot hold at an optimal solution.
We will show that this holds if 
\begin{equation}
\label{eq:norefw2}
w^{{\tt ref}, \eta}< \wrefp,\ \ w^{{\tt ref}, \eta}< z_j \enspace j \in \{1,\ldots,n\},
\end{equation}
which is certainly true for $h$ small enough, recalling that $w^{{\tt ref}, \eta}$ converges to 0 as $h\rightarrow 0$.

Recalling that we are addressing the case $f_{j-1}>0, f_j<0$,
we first consider the subcase $j=2$. In such case we have $w_{j-1}=\winit=z_1$ and 
\begin{align*}
f_1&=\frac{w_2-z_1}{h}+\gamma z_1 -g(\sin(\alpha_1) + c\cos(\alpha_1))=\\
&=\frac{w^{{\tt ref}, \eta}-z_1}{h}+\gamma z_1 -g(\sin(\alpha_1) + c\cos(\alpha_1)),
\end{align*}
which is lower than $-g\mu\cos(\alpha_1)$ for $h$ small enough, thus violating the minimum force constraint $f_1\geq -g\mu\cos(\alpha_1)$ in~\eqref{eq:probfix}.
\noindent
If $j>2$ we also consider force $f_{j-2}$. We have three possibilities:
\begin{itemize}
\item If $f_{j-2}<0$, in view of case {\bf b)} of the current proposition with $f_{j-2}<0, f_{j-1}>0$, we must have $w_{j-1}=z_{j-1}$.
Therefore:
\begin{align*}
f_{j-1}&
=\frac{w_j-w_{j-1}}{h}+\gamma w_{j-1} -g(\sin(\alpha_{j-1}) + c\cos(\alpha_{j-1}))= \\
& =\frac{w^{{\tt ref}, \eta}-z_{j-1}}{h}+\gamma z_{j-1} -g(\sin(\alpha_{j-1}) + c\cos(\alpha_{j-1})),
\end{align*}
which is lower than $-g\mu\cos(\alpha_{j-1})$ for $h$ small enough, thus violating the minimum force constraint $f_{j-1}\geq -g\mu\cos(\alpha_{j-1})$ in~\eqref{eq:probfix}.
\item If $f_{j-2}>0$, in view of case {\bf d)} of the current proposition with $f_{j-2}>0, f_{j-1}>0$, we have that $w_{j-1}=z_{j-1}$ or $w_{j-1}=y_{j-1}\geq \wrefp$ or $w_{j-1}=\wrefp$.
Again, for $h$ small enough we have that $f_{j-1}$ violates the minimum force constraint in~\eqref{eq:probfix}.
\item If $f_{j-2}=0$, in view of case {\bf h)} of the current proposition with $f_{j-2}=0, f_{j-1}>0$, we have that $w_{j-1}=z_{j-1}$ or $w_{j-1}\geq \wrefp$. In both cases we have the same violation as before of the minimum force
constraint $f_{j-1}\geq -g\mu\cos(\alpha_{j-1})$.
\end{itemize}
\end{proof}

\section{A dynamic programming approach}
\label{sec:dynprog}
In this section, we propose a possible DP approach which explores a subset of the feasible solutions of Problem~\eqref{eq:newconvrel}, and we will prove that 
the optimal solution of~\eqref{eq:newconvrel} is close to such subset of feasible solutions.

\subsection{State definition} 
The results of Proposition~\ref{prop:combposs} allow us to
restrict the attention to four possible speeds at each step $j \in \{1,\ldots,n\}$:
\begin{itemize}
\item[i)] the minimum speed $y_j$;
\item[ii)] the maximum speed $z_j$;
\item[iii)] the reference speed $\wrefp$;
\item[iv)] the reference speed $\wrefm$.
\end{itemize}
Note that these are not the only four speeds that we consider.
In fact, these speeds are connected to each other by sections of speed profiles associated with null-force, as we will explain below.
We denote by $W_j=\{y_j,z_j,\wrefp,\wrefm\}$ the sets of four speeds at step $j$. 
In the proposed DP approach we employ the following set of states:
\[
\Omega=\left\{(j,w_j)\ \mid\ j \in \{1,\ldots,n\},\ w_j\in W_j\right\}.
\]
Note that each state has two components, the first one being a step $j$, while the second component is one of the four speeds in $W_j$.
Note that $|\Omega|=4n$. There is a single initial state $(1,\winit)$, and a single final state $(n,\wfin)$.

\subsection{Moves}
For a given speed $w_j$ at step $j$, we define the set of feasible speeds at step $j+1$, with respect to the constraints in~\eqref{eq:probfix}, as follows:
\begin{align*}
\Gamma(w_j)=&\left\{w_{j+1}\in [y_{j+1},z_{j+1}]\ \mid\ f_j\geq -g\mu\cos(\alpha_j), \phantom{\frac11}\right.\\
& \left. f_{j} \leq \min\left\{g\mu\cos(\alpha_j),\frac{\Pmax}{M\sqrt{2w_j}}\right\}\right\}.
\end{align*} 
In other words, $\Gamma(w_j)$ is the set of all values which can be assigned to $w_{j+1}$ in such a way that $w_{j+1}$ fulfills the lower and upper limit constraints for such variable,
and force $f_j$, defined, as usual, as follows:
\[
f_j=\frac{1}{h} (w_{j+1}-w_j) + \gamma w_j + g (\sin \alpha_j + c\cos(\alpha_j)),
\]
fulfills the minimum and maximum force constraints at step $j$, and the power constraint at step $j$. 
From each state we define two distinct types of moves: {\em single-step} moves and {\em null-force} moves. 
Both types of moves generate a state where the second component is one of the previously defined four speeds, starting from a state where, again, the second component is one of such four speeds.

{\bf Single-step moves:} For a given state $(j,w_j)$, $w_j\in W_j$, single-step moves are those leading to the set of states:
\[
S_1(j,w_j)=\left\{(j+1,w_{j+1})\ \mid\ w_{j+1}\in W_{j+1}\cap \Gamma(w_j)\right\}.
\]
Note that $|S_1(j,w_j)|\leq 4$.

{\bf Null-force moves:} These are moves starting at some step $j$ and obtained through a sequence of steps where the force is null.
First, we define the {\em null-force speed curve} starting at step $j$ with speed $w_j$, that is, the speed values $w_j,\ldots,w_n$ which are obtained by starting at speed $w_j$ and setting $f_j=\ldots=f_{n-1}=0$:
\begin{equation}
\label{eq:nullforcurve}
\begin{aligned}
w_j^{j,0}(w_j)=&\ w_j,\\
w_{j+i}^{j,0}(w_j)=&\ (1-h\gamma)w_{j+i-1}^{j,0}(w_j) -h g (\sin(\alpha_{j+i-1}) + c\cos(\alpha_{j+i-1})),\\
&\text{for }i \in \{1,\ldots,n-j\}.
\end{aligned}
\end{equation}
We define a set of distinct states which can be reached from some state $(j,w_j)$, $w_j\in W_j$, by a null-force move. We first compute index $\eta(j,w_j)$ as the maximum number of steps after $j$ where the null-force speed curve does not violate the lower or upper bound limits for the speeds:
\begin{align*}
\eta(w_j,j)=\max &\left\{ i\in \{0,\ldots,n-j-1\}\ \mid w_{j+i}^{j,0}(w_j) \in [y_{j+i},z_{j+i}]\right\}.
\end{align*}
Next, we define the sets of indices:
\begin{equation}
\label{eq:defI}
\begin{aligned}
& I_1(w_j,j)=\left\{i\in\{1,\ldots,\eta(w_j,j)\} \mid z_{j+i+1}\in \Gamma(w^{j,0}_{j+i}(w_j))\right\} \\
& I_2(w_j,j)=\left\{i\in\{1,\ldots,\eta(w_j,j)\} \mid y_{j+i+1}\in \Gamma(w^{j,0}_{j+i}(w_j))\right\} \\
& I_3(w_j,j)=\left\{i\in\{1,\ldots,\eta(w_j,j)\} \mid \wrefp\in \Gamma(w^{j,0}_{j+i}(w_j))\right\} \\
& I_4(w_j,j)=\left\{i\in\{1,\ldots,\eta(w_j,j)\} \mid \wrefm\in \Gamma(w^{j,0}_{j+i}(w_j))\right\},
\end{aligned}
\end{equation}
that is, all steps after $j$ such that, from the corresponding point along the null-force curve, the maximum speed profile $z$ (index set $I_1$), or the minimum speed profile (index set $I_2$), or the reference speed $\wrefp$ (index set $I_3$),
or the reference speed $\wrefm$ (index set $I_4$), can be reached by a feasible one-step move.
Note that some of these sets (but not all) might be empty.
For each $i\in I_1(w_j,j)$ we define the new state $(j+i+1,z_{j+i+1})$ as a result of a null-force move. Similarly, for each $i\in I_2(w_j,j), I_3(w_j,j), I_4(w_j,j)$, we define a new state, as previously described, where the second component of the new state becomes $y_{j+i+1}$, $\wrefp$, and $\wrefm$, respectively.
If we set:
\[
\tilde{w}_{j+i+1}=
\begin{cases}
z_{j+i+1} & \mbox{if } i\in I_1(w_j,j) \\
y_{j+i+1} & \mbox{if } i\in I_2(w_j,j) \\
\wrefp & \mbox{if } i\in I_3(w_j,j) \\
\wrefm & \mbox{if } i\in I_4(w_j,j),
\end{cases}
\]
then, the set of possible states which can be reached by null-force moves starting from $(j,w_j)$ is:
\[
S_{f=0}(j,w_j)=\bigcup_{r=1}^4 \left\{(j+i+1,\tilde{w}_{j+i+1})\ \mid\ i\in I_r(w_j,j)\right\}.
\]
Note that in each of these states the second component is a speed belonging to the corresponding set of four speeds $W_{j+i+1}$.
The following proposition shows that the sequence of speeds generated by a null-force move is feasible with respect to all the corresponding constraints in~\eqref{eq:probfix}.
\begin{prop}
\label{prop:feasseq}
The sequence of speeds $w_{j}^{j,0}(w_j)$,\ldots, $w_{j+i}^{j,0}(w_j)$, $\tilde{w}_{j+i+1}$, for $i\in I_r(w_j,j)$, $r \in \{1,\ldots,4\}$, is feasible with respect to all the corresponding constraints in~\eqref{eq:probfix}.
\end{prop}
\begin{proof}
We notice that, by definition of $\eta(w_j,j)$, $w_{j+i}^{j,0}(w_j)\in [y_{j+i},z_{j+i}]$ for $i \in \{0,\ldots,\eta(w_j,j)\}$. Moreover, all the forces $f_{j+r}$, $r \in \{0,\ldots, i-1\}$ are equal to $0$, so that the force and power constraints are fulfilled. Finally, by definition of each set $I_r(w_j,j)$, $r \in \{1,\ldots,4\}$, we have that 
$\tilde{w}_{j+i+1}\in \Gamma(w^{j,0}_{j+i}(w_j))$, so that
the force $f_{j+i}$ fullfills the corresponding power and force constraints.
\end{proof}

\subsection{Dynamic programming approach}
Now we are ready to introduce the DP approach which returns an (approximate) solution of Problem~\eqref{eq:probfix}.
Algorithm~\ref{alg:dynprog} returns the objective function value of the solution, while Algorithm~\ref{alg:solbuild} allows building the (approximate) optimal solution starting from the output of
Algorithm~\ref{alg:dynprog}.

\begin{algorithm}
\caption{\label{alg:dynprog} Dynamic programming procedure}
\begin{algorithmic}[1]
\Procedure{DynProg}{$y,z,\wrefp,\wrefm$}
\State Initialize the set of states ${\mathcal{Q}}=\{(1,\winit)\}$ \label{lin:1dyn}\;
\State Set $V(1,\winit)=0$, $P(1,\winit)=-$ \label{lin:1bisdyn} \;
\While{${\mathcal{Q}}\neq \varnothing$}
\State Extract a state $(j,w_j)$ from ${\mathcal{Q}}$ with lowest first component \label{lin:2dyn}\;
\State Set ${\mathcal{Q}}={\mathcal{Q}}\setminus \{(j,w_j)\}$\label{lin:3dyn}\;
\For{$(j+1,w_{j+1}) \in S_1(j,w_j)$}
\State Let $J(j+1,w_{j+1})=V(j,w_j)+\lambda M \eta \gamma w_j+\frac{1}{\sqrt{2w_j}}+(1-\eta)\lambda M \max\{f_j,0\}$\label{lin:4dyn}\;
\If{$(j+1,w_{j+1})\not\in {\mathcal{Q}}$}
\State Add $(j+1,w_{j+1})$ to ${\mathcal{Q}}$ \label{lin:5dyn}\;
\State Set $V(j+1,w_{j+1})=J(j+1,w_{j+1})$, $P(j+1,w_{j+1})=(j,w_j)$ \label{lin:7dyn}\;
\ElsIf{$J(j+1,w_{j+1})< V(j+1,w_{j+1})$}
\State Set $ V(j+1,w_{j+1})=J(j+1,w_{j+1})$, $P(j+1,w_{j+1})=(j,w_j)$ \label{lin:9dyn}\;
\EndIf
\EndFor
\For{$(j+i,w_{j+i}) \in S_{f=0}(j,w_j)$}
\State Let $J(j+i,w_{j+i})=V(j,w_j)+\left[\sum_{r=0}^{i-1} \lambda M \eta \gamma w_{j+r}+\frac{1}{\sqrt{2w_{j+r}}}\right]+(1-\eta)\lambda M \max\{f_{j+i-1},0\}$\label{lin:10dyn}\;
\If{$(j+i,w_{j+i})\not\in {\mathcal{Q}}$}
\State Add $(j+i,w_{j+i})$ to ${\mathcal{Q}}$ \label{lin:11dyn}\;
\State Set $V(j+i,w_{j+i})=J(j+i,w_{j+i})$, $P(j+i,w_{j+i})=(j,w_j)$\label{lin:12bisdyn}\;
\ElsIf{$J(j+i,w_{j+i}) < V(j+i,w_{j+i})$}
\State Set $ V(j+i,w_{j+i})=J(j+i,w_{j+i})$, $P(j+i,w_{j+i})=(j,w_j)$\label{lin:15dyn}\;
\EndIf
\EndFor
\EndWhile
\State \Return{$V(n,\wfin), P$}\;
\EndProcedure
\end{algorithmic}
\end{algorithm}

\noindent
In line~\ref{lin:1dyn} we initialize a queue of states ${\mathcal{Q}}$ with the single initial state $(1,\winit)$. 
In line~\ref{lin:1bisdyn} we set the value of the initial state equal to $0$ and leave the predecessor $P$ of such state undefined.
The function value $V$ of a state $(j,w_j)$ is the best possible cumulative objective function value observed up to step $j-1$ (i.e., the first $j-1$ terms of function $F$) to reach
 state $(j,w_j)$ from initial state $(1,\winit)$, while $P$ associates to a state $(j,w_j)$, state $(j-r,w_{j-r})$ for some $r>0$, from which a transition to $(j,w_j)$ has been performed, delivering the current best value
 $V(j,w_j)$ for that state. Then, we enter a {\tt While} loop that repeats until queue ${\mathcal{Q}}$ is empty.
 In line~\ref{lin:2dyn} we extract a state in the queue with the lowest possible first component, and in line~\ref{lin:3dyn} we remove such state from the queue.
 Then, in the first {\tt For} cycle, for each state $(j+1,w_{j+1}) \in S_1(j,w_j)$ we compute the cumulative objective function value $J(j+1,w_{j+1})$
 given by the value $V(j,w_j)$ of state $(j,w_j)$ and the contribution to the objective function value of step $j$ (line~\ref{lin:4dyn}). If state $(j+1,w_{j+1})$ is not already present in the queue, we add it
 to ${\mathcal{Q}}$ and assign value $J(j+1,w_{j+1})$ and predecessor $(j,w_j)$ to it (lines~\ref{lin:5dyn}--\ref{lin:7dyn}).
 Instead, if $(j+1,w_{j+1})$ already belongs to ${\mathcal{Q}}$, we compare $J(j+1,w_{j+1})$ with $V(j+1,w_{j+1})$ and if the former is lower than the latter, we update 
 $V(j+1,w_{j+1})$ and $P(j+1,w_{j+1})$ (line~\ref{lin:9dyn}).
 In the second {\tt For} cycle we repeat the same operations of the first {\tt For} cycle with the states in $S_{f=0}(j,w_j)$.
 The algorithm returns $V(n,\wfin)$, which is the (approximate) optimal value of Problem~\eqref{eq:probfix} and  function $P$, which is needed to build the (approximate) optimal solution.
 More precisely, the approximate optimal solution is built through Algorithm~\ref{alg:solbuild}.
 In line~\ref{lin:1build} we initialize index $j$ with $n$ and speed $\bar{w}$ with $\wfin$. We also set the $n$-th component of the approximate solution $w^\star$ equal to $\wfin$.
 Then, we enter a {\tt While} loop from which we exit as soon as we reach step $j=1$. In line~\ref{lin:2build} we define the predecessor $(j-r,w_{j-r})$ of the current state $(j,\bar{w})$.
 In case the first component of the predecessor is $j-1$ (i.e., we reached state $(j,\bar{w})$ through a single-step move), then in line~\ref{lin:3build} we set the $(j-1)$-th component of $w^\star$ equal to the second component
 of the predecessor. Moreover, we set $j:=j-1$ and $\bar{w}=w_{j-1}$ (line~\ref{lin:3build}), that is, we set the current state equal to the predecessor $(j-1,w_{j-1})$.
 Instead, if the first component of the predecessor is $j-r$ with $r>1$, we reached state $(j,w_j)$ with a null-force move starting at $w_{j-r}$. We set all speeds
 $w^\star_{j-r+i}$, $i \in \{1,\ldots,r-1\}$, by starting at $w^\star_{j-r}=w_{j-r}$ (see line~\ref{lin:5build}), and moving along the null-force curve 
 (see line~\ref{lin:6build} in the body of the {\tt For} cycle).
 Finally, in line~\ref{lin:7build} we update the current state by setting it equal to the predecessor $(j-r,w_{j-r})$. 

\begin{algorithm}
\caption{\label{alg:solbuild} Procedure to build the (approximate) optimal solution}
\begin{algorithmic}[1]
\Procedure{BuildSolution}{$P$}
\State Set $j=n$, $\bar{w}=\wfin$, $w_n^\star=\wfin$ \label{lin:1build}\;
\While{$j\neq 1$}\;
\State Let $(j-r,w_{j-r})=P(j,\bar{w})$ \label{lin:2build}\;
\If{$r=1$}
\State Set $w^\star_{j-1}=w_{j-1}$\label{lin:3build}\;
\State Set $\bar{w}=w_{j-1}$, $j=j-1$ \label{lin:4build}\;
\Else
\State $w_{j-r}^\star=w_{j-r}$ \label{lin:5build}\;
\For{$i \in \{1,\ldots,r-1\}$}
\State Set $w_{j-r+i}^\star=(1-h\gamma)w_{j-r+i-1}^\star -hg(\sin(\alpha_{j-r+i-1}) + c\cos(\alpha_{j-r+i-1}))$ \label{lin:6build}\;
\EndFor
\State Set $\bar{w}=w_{j-r}$, $j=j-r$ \label{lin:7build}\;
\EndIf
\EndWhile
\State \Return{$w^\star$}\;
\EndProcedure
\end{algorithmic}
\end{algorithm}

The following proposition establishes the time-complexity of Algorithms~\ref{alg:dynprog} and~\ref{alg:solbuild}.
\begin{prop}
\label{prop:compl}
The time-complexity of Algorithm~\ref{alg:dynprog} is $O\left(n^2\right)$.
The time-complexity of Algorithm~\ref{alg:solbuild} is $O\left(n\right)$.
\end{prop}
\begin{proof}
First of all, we observe that the maximum number of states is $4n$. For each state,
we perform at most four one-step moves which require $O(1)$ operations. Moreover, for each state we perform a null-force move, which requires $O(n)$ operations. 
Therefore, the overall number of operations is
$O(n^2)$.
Concerning Algorithm~\ref{alg:solbuild}, it moves backward from step $n$ to step $1$, and we notice that for each step $j$ the number of operations is $O(1)$, so that the overall number of operations
is $O(n)$.
\end{proof}

\subsection{Distance between the optimal solution and the solutions explored by the DP approach}
While the proposed DP approach is quite efficient, it does not deliver the exact optimal solution of Problem~\eqref{eq:probfix}, since
at each step $j$ it only considers the four speeds in $W_j$, while speeds outside sets $W_j$ can only be visited along null-force moves. 
However, we can prove that the difference between the speeds at optimal solutions of Problem~\eqref{eq:probfix} and at least one of those explored by the DP approach is bounded from above by 
a quantity proportional to $h$.
Moreover, also the values of the two objective functions are close.
This is stated in the following final result, whose proof is presented in Appendices~\ref{sec:prelim_res}--\ref{sec:final_proof}.
\begin{thrm}
\label{thrm:main}
Let $(w^\star,f^\star)$ be the optimal solution of Problem~\eqref{eq:probfix}. Then, there exists $(\bar{w},\bar{f})\in \bar{W}$ such that
\[
|\bar{w}_j-w_j^\star|=O(h)\qquad j\in\{1,\ldots,n\}.
\]
Moreover, let
\[
{\mathcal{J}}_{\bar w}=\{j\ \mid\ \bar{w}_j\neq w^\star_j\},\qquad {\mathcal{J}}_{\bar f}=\{j\ \mid\ \bar{f}_j\neq f^\star_j\},
\]
be the set of indices where $\bar{w}$ ($\bar{f}$, respectively) is different from $w^\star$ ($f^\star$, respectively).
If $w_j^\star>0$ for all $j\in {\mathcal{J}}_{\bar w}$,
then
\[
\left| hF(\bar{w},\bar{f})-hF(w^\star,f^\star)\right|\leq O(h)+h \eta\gamma \lambda g M \mu |{\mathcal{J}}_f|.
\]
\end{thrm}

Note that a detailed discussion on how to build solution $(\bar{w},\bar{f})$ of Theorem~\ref{thrm:main} is provided in Appendix~\ref{sec:solclose}.

\section{Computational experiments}
\label{sec:compexp}
Algorithm~\ref{alg:dynprog} was implemented in MATLAB~\cite{MATLAB2024} and converted to C using the MATLAB Coder R2024b. The tests were performed on an Apple M2 Pro system with 10 cores (6 performance cores and 4 efficiency cores) and \SI{16}{\giga\byte} of RAM.
When mentioned, the exact convex relaxation~\eqref{eq:newconvrel} of Problem~\eqref{eq:probfix} was solved using the commercial solver MOSEK (version 10.1.31)~\cite{mosek}, interfaced with MATLAB via YALMIP~\cite{lofberg2004yalmip}. In this case, only MOSEK solution times were considered, while the YALMIP overhead was neglected.

\begin{figure}[h!]
	\centering
	\includegraphics[width=.6\columnwidth]{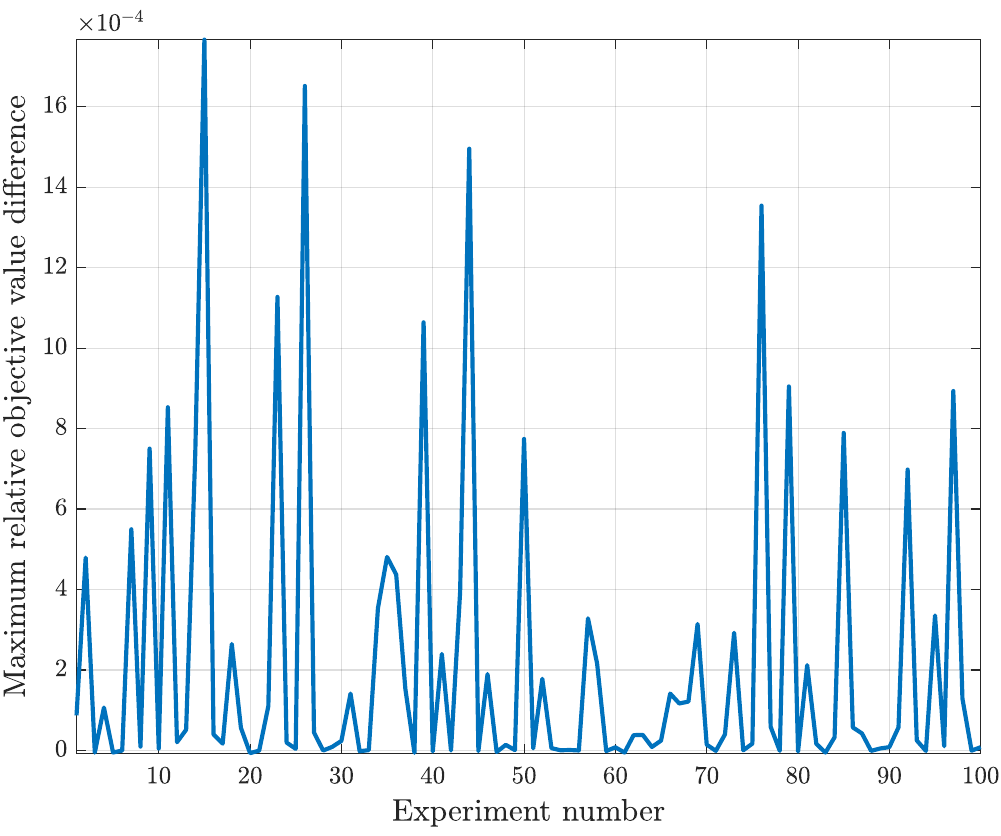}
	\caption{Maximum relative objective function value differences.}
	\label{fig:objDRel}
\end{figure}

\begin{figure}[h!]
	\centering
	\includegraphics[width=.6\columnwidth]{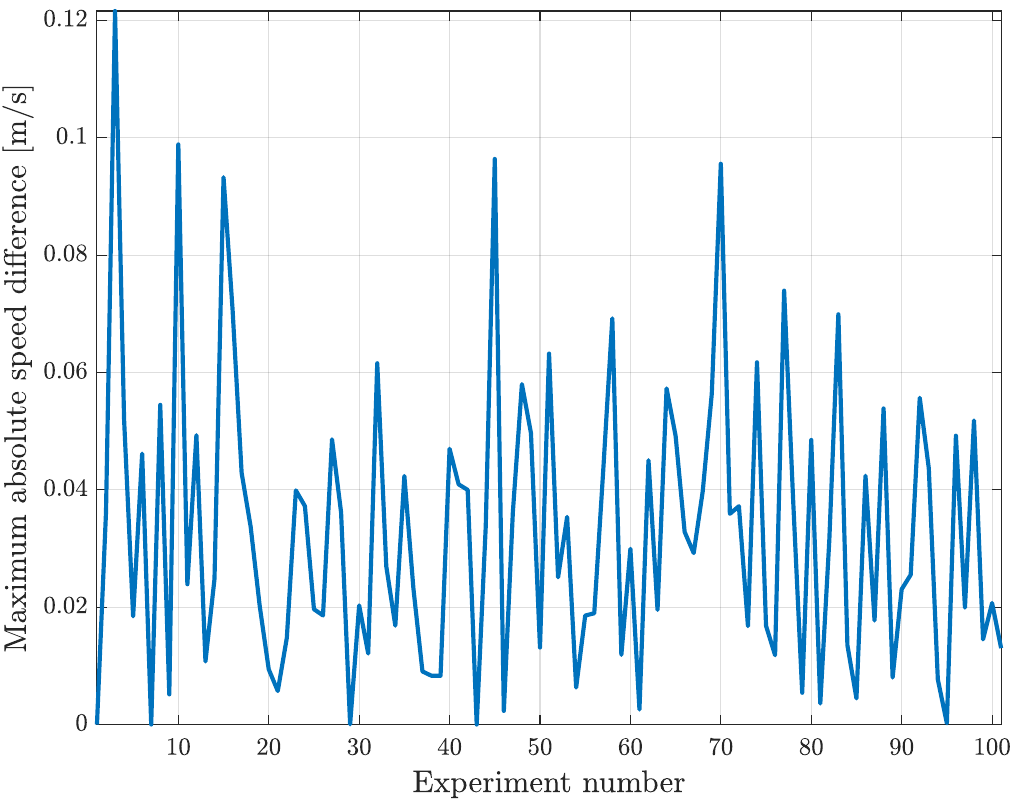}
	\caption{Maximum absolute speed differences.}
	\label{fig:vDiff}
\end{figure}

In all experiments, the planning is addressed considering a Fiat 500e full electric, whose parameters are reported in Table~\ref{tab:vehicle_spec}.
Moreover, we considered $\mu = 0.7$, representing the static friction coefficient between tires and dry asphalt
\begin{table}[h!]
\centering
\begin{tabular}{|c|c|}
\hline
$M$			& \SI{1365}{\kilogram}					\\
\hline
$\Pmax$		& 118\ metric hp $\approx$ \SI{87}{\kilo\watt}	\\
\hline
$\vmax$		& \SI{150}{\kilo\meter\per\hour} 			\\
\hline
$\wmax$		& \SI{868}{\meter^2\per\second^2}			\\
\hline
\end{tabular}
\quad
\begin{tabular}{|c|c|}
\hline
$\eta$		& $70\%$ \\
\hline
$c$		& 0.007 \\
\hline
$\Gamma$	& \SI{0.399}{\kilogram\per\meter} \\
\hline
$\gamma$	& \SI{0.2923e-3}{\per\meter} \\
\hline
\end{tabular}
\caption{Fiat 500e specifications.}
\label{tab:vehicle_spec}
\end{table}
First, we tested the efficacy of Algorithm~\ref{alg:dynprog} in terms of both computational times and accuracy. We fixed $n = 2000$ and $h =  \SI{0.2}{\meter}$ and generated 100 random instances of Problem~\eqref{eq:probfix} with piecewise constant maximum speed and slope, and randomly selected the initial and final speeds. We compared the computational performance of Algorithm~\ref{alg:dynprog} with that of the aforementioned commercial solver.
\begin{figure}[h!]
	\centering
	\includegraphics[width=.6\columnwidth]{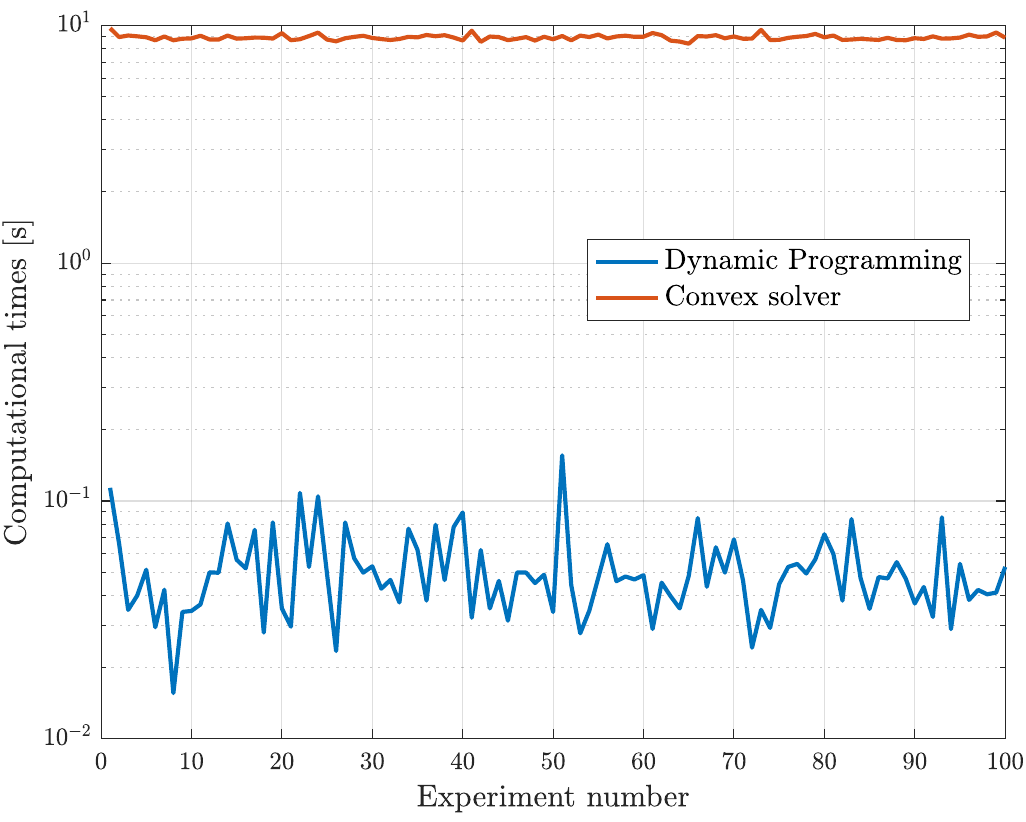}
	\caption{Computational times of the DP and the commercial solver.}
	\label{fig:times}
\end{figure}
Figure~\ref{fig:times} shows that Algorithm~\ref{alg:dynprog} is approximately 3 orders of magnitude faster than the commercial solver. This significant improvement in computational time comes at the expense of a slight degradation of the objective function value, whose maximum relative differences are presented in Figure~\ref{fig:objDRel}.
As we can see, the increase in the objective function values when employing Algorithm~\ref{alg:dynprog} are negligible. From the perspective of the speed profiles, the maximum absolute differences between the profiles obtained with Algorithm~\ref{alg:dynprog} and those computed by the commercial solver are moderate and acceptable (see Figure~\ref{fig:vDiff}).

\begin{figure}[h!]
	\centering
	\includegraphics[width=.6\columnwidth]{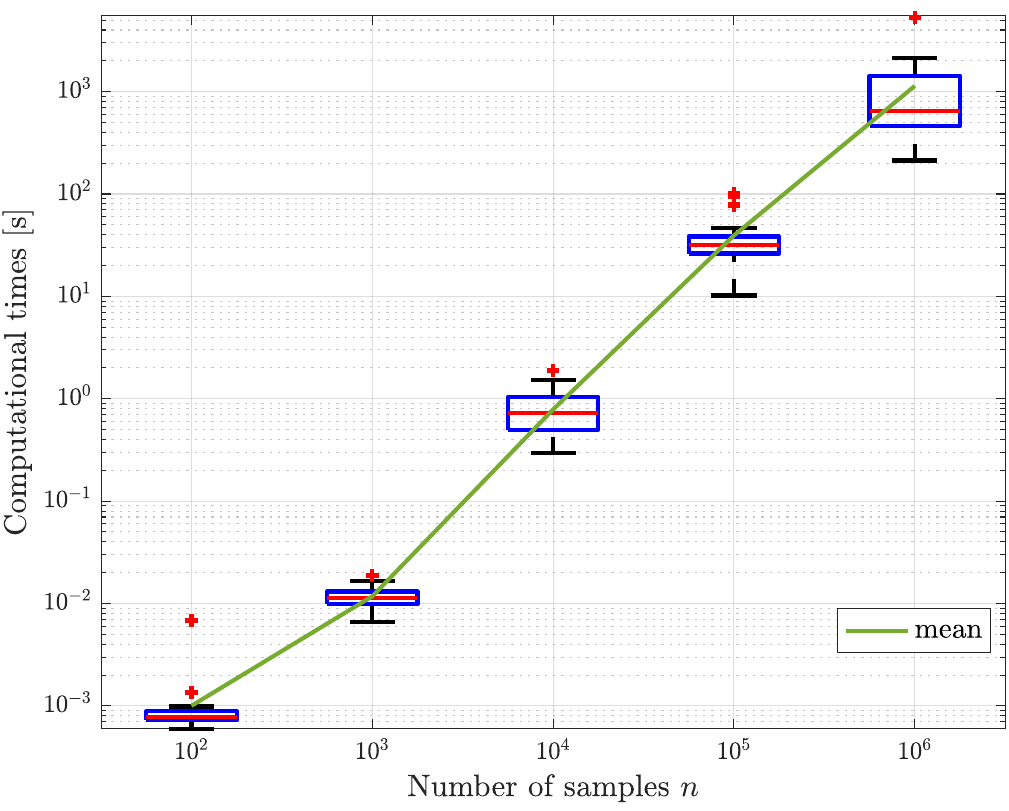}
	\caption{Computational times of the DP for different values of $n$.}
	\label{fig:nRange}
\end{figure}

We also evaluated the computational performance of Algorithm~\ref{alg:dynprog} for different values of $n$. We considered values of $n$ in range $[10^2, 10^6]$, using a logarithmic scale with 5 samples. For each value of $n$, we generated 30 random instances of Problem~\eqref{eq:probfix} as in the previous set of experiments. As shown in Figure~\ref{fig:nRange}, the performance of Algorithm~\ref{alg:dynprog} is compatible with real-time applications for values of $n$ up to the order of $10^3$.

\begin{figure}[h!]
	\centering
	\begin{subfigure}{0.49\columnwidth}
	\centering
	\includegraphics[width=\linewidth]{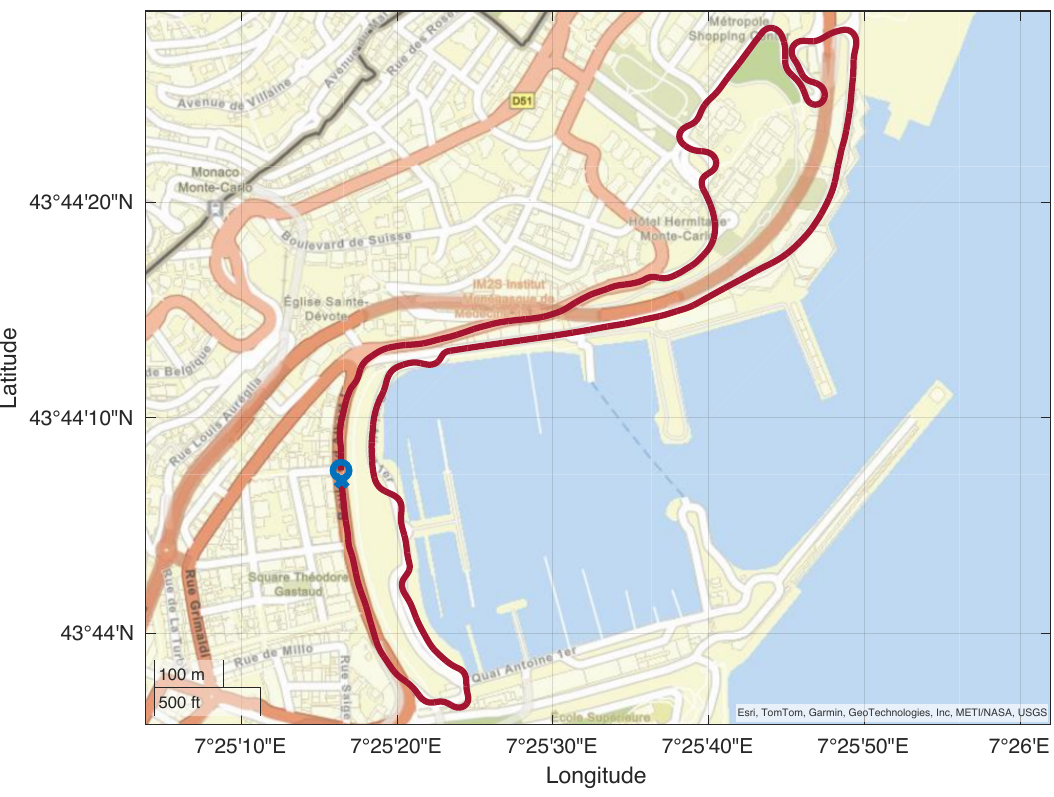}
	\caption{Monaco Grand Prix circuit.}
	\label{fig:MonacoCircuit}
    \end{subfigure}
    \hfill
    \begin{subfigure}{0.49\columnwidth}
        \centering
        \includegraphics[width=\linewidth]{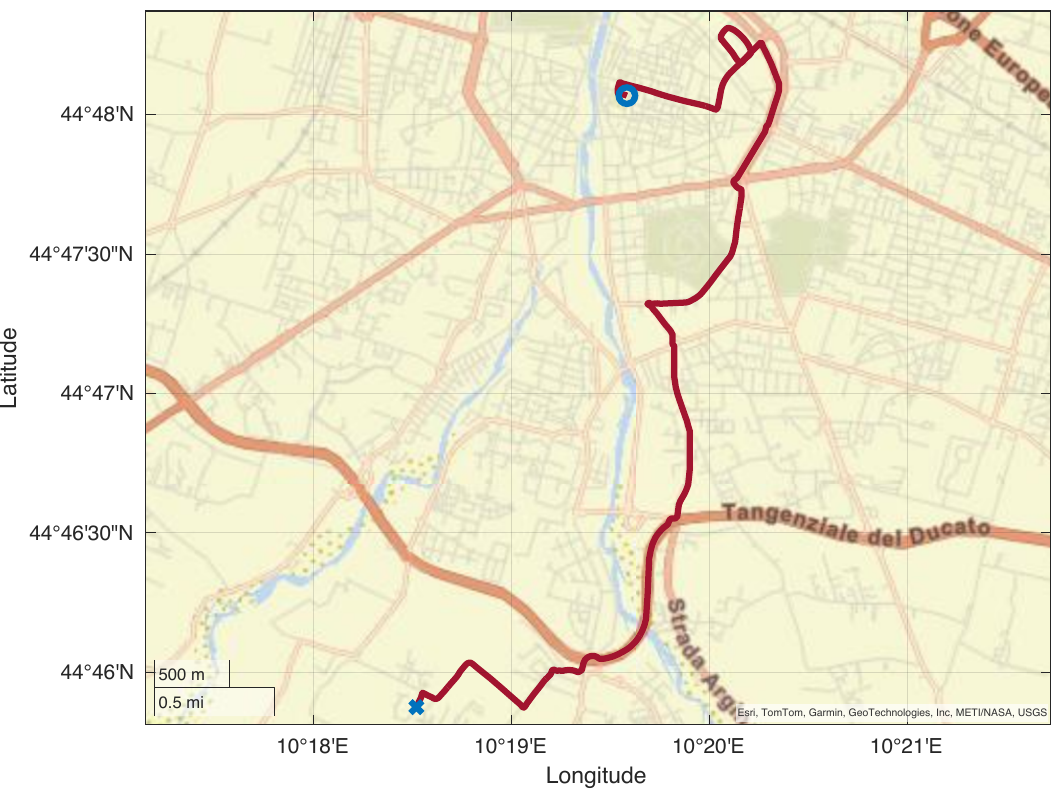}
        \caption{Urban path in Parma.}
        \label{fig:UniPR}
    \end{subfigure}
	\begin{subfigure}{0.49\columnwidth}
	\centering
	\includegraphics[width=\linewidth]{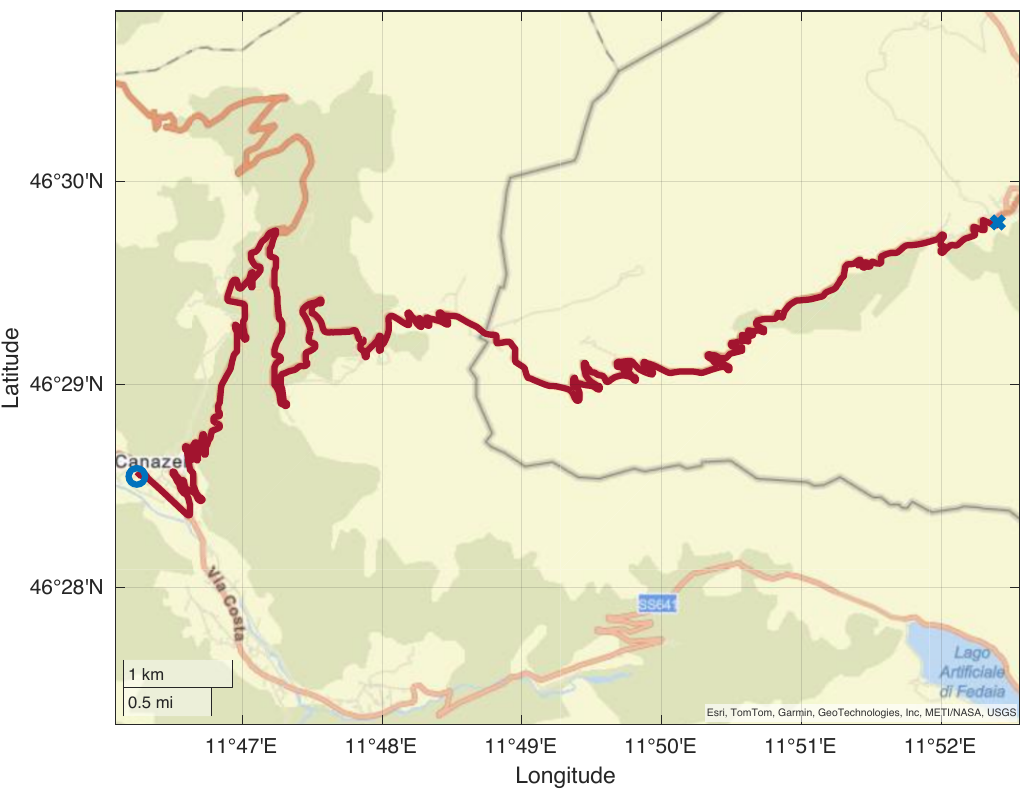}
	\caption{Passo Pordoi.}
	\label{fig:Arabba-Canazei}
    \end{subfigure}
    \hfill
    \begin{subfigure}{0.49\columnwidth}
        \centering
        \includegraphics[width=\linewidth]{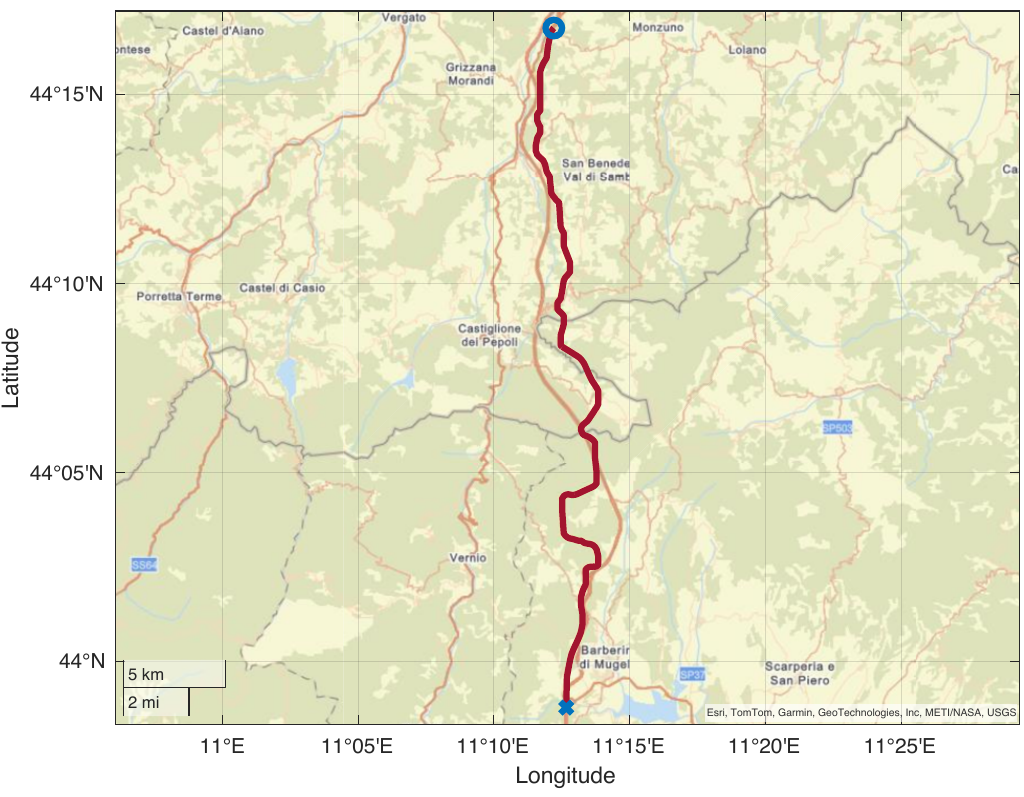}
        \caption{Section of Autostrada A1.}
        \label{fig:Rioveggio-Barberino}
    \end{subfigure}
    \caption{Example paths with circles and crosses denoting starting and ending points, respectively.}
    \label{fig:paths}
\end{figure}

\begin{figure}[h!]
	\centering
	\begin{subfigure}{0.49\columnwidth}
	\centering
	\includegraphics[width=\linewidth]{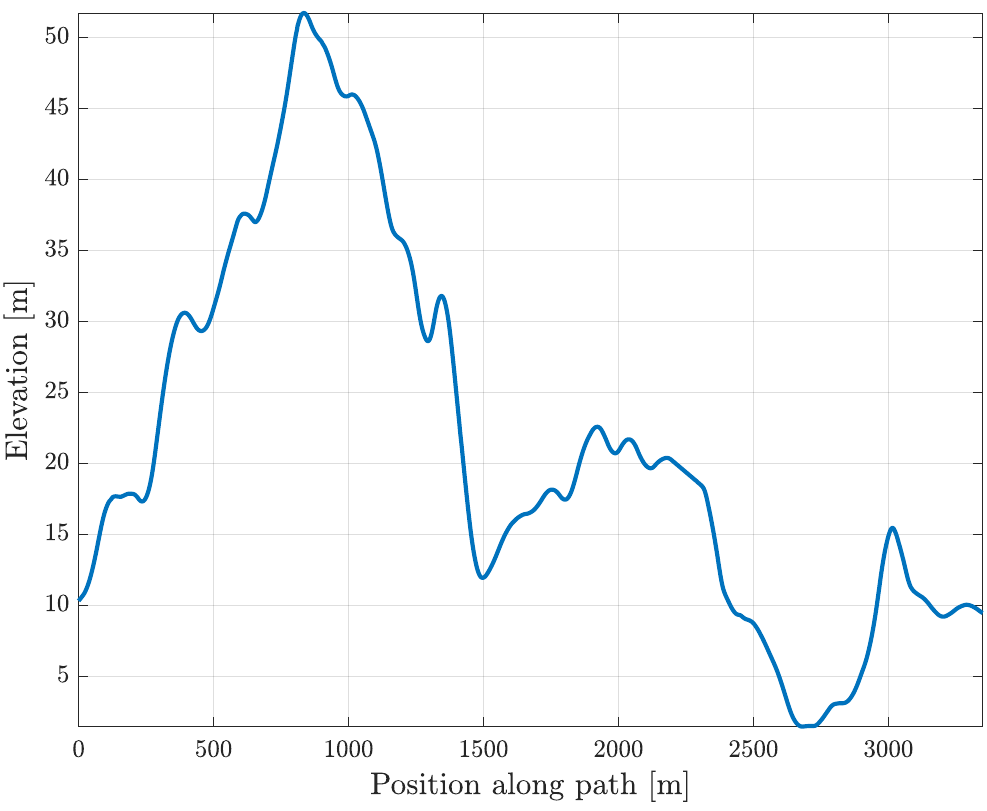}
	\caption{Monaco Grand Prix circuit.}
	\label{fig:MonacoElevation}
    \end{subfigure}
    \hfill
    \begin{subfigure}{0.49\columnwidth}
        \centering
        \includegraphics[width=\linewidth]{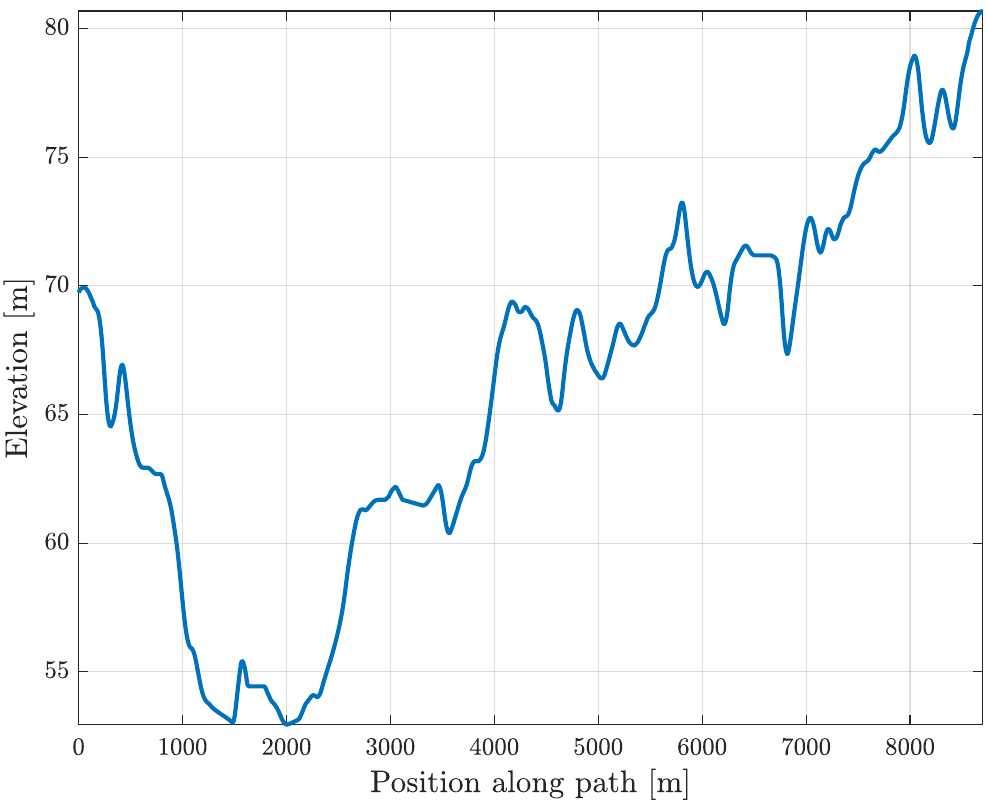}
        \caption{Urban path in Parma.}
        \label{fig:UniPRElevation}
    \end{subfigure}
	\begin{subfigure}{0.49\columnwidth}
	\centering
	\includegraphics[width=\linewidth]{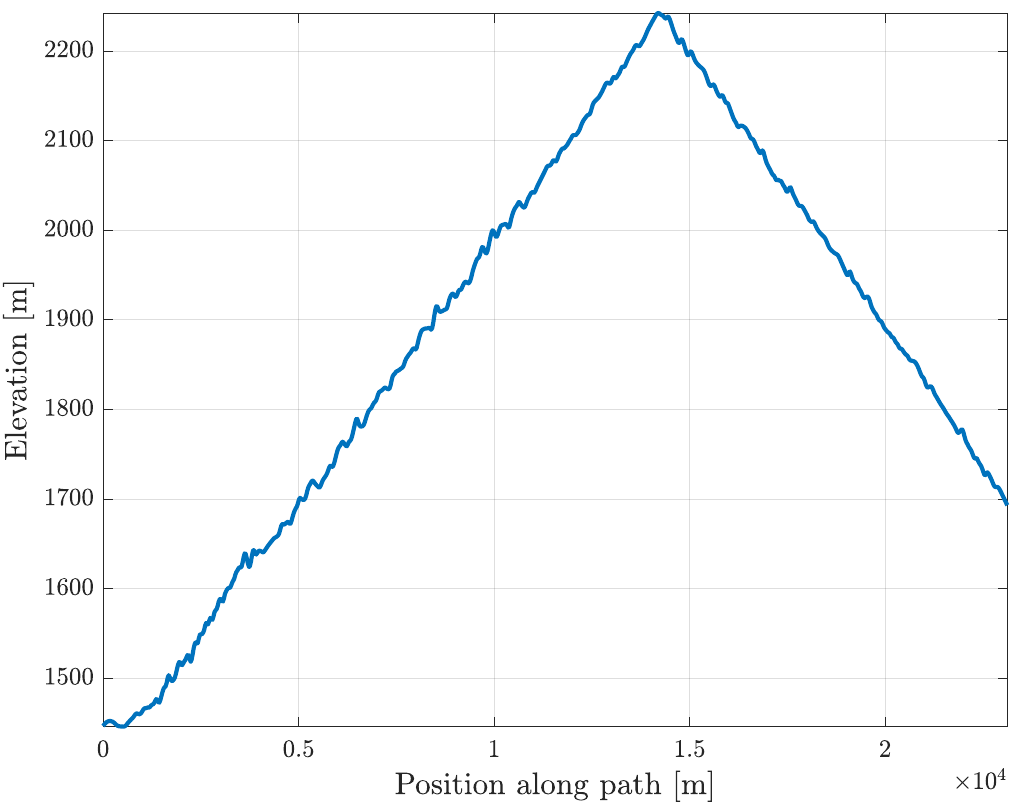}
	\caption{Passo Pordoi.}
	\label{fig:Arabba-CanazeiElevation}
    \end{subfigure}
    \hfill
    \begin{subfigure}{0.49\columnwidth}
        \centering
        \includegraphics[width=\linewidth]{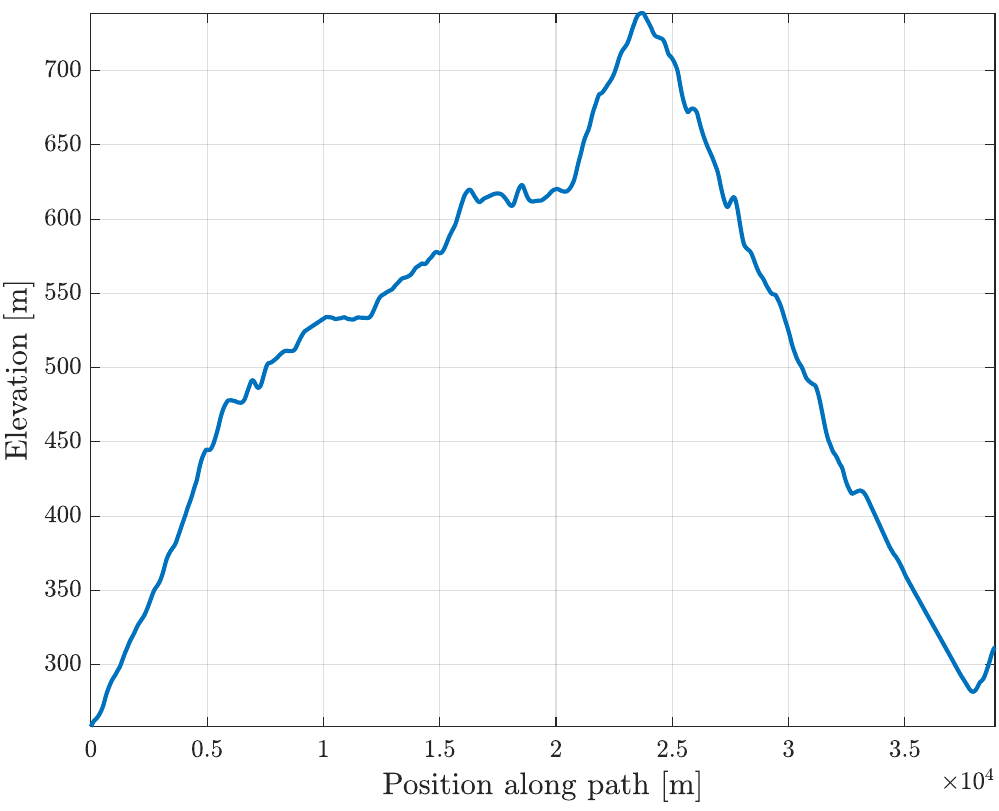}
        \caption{Section of Autostrada A1.}
        \label{fig:Rioveggio-BarberinoElevation}
    \end{subfigure}    
    \caption{Elevation profiles.}
    \label{fig:elevations}
\end{figure}

Finally, in order to test the robustness of the proposed algorithm across a variety of driving conditions, we evaluate it on four distinct scenarios (see Figure~\ref{fig:paths} for the layouts and Figure~\ref{fig:elevations} for the corresponding elevation profiles, respectively):
the Monaco Grand Prix circuit, a short temporary street circuit characterized by tight turns and frequent direction changes;
an urban path in Parma, Italy, passing by several University Departments, with many tight corners and low speeds;
a steep and winding mountain pass (Passo Pordoi, Italy) with significant elevation and many hairpin turns;
a hilly highway section (Autostrada A1, Italy, from Rioveggio to Barberino di Mugello) with a longer distance, higher speeds, with uphill and downhill parts, and more straight segments.
In all four scenarios the vehicle starts and ends its travel at zero speed, that is $w_1 = w_n = 0$ \unit{\meter\squared\per\second\squared}.

We sampled the values of $\lambda$ in range $[0, 1]$ using a logarithmic scale with 45 samples. For each value, we solved optimization Problem~\eqref{eq:probfix} by means of Algorithm~\ref{alg:dynprog} and computed the corresponding optimal values of travel time and specific energy consumption (i.e., the amount of energy consumed by the vehicle, normalized by its mass).

\begin{figure}[h!]
	\centering
	\begin{subfigure}{0.49\columnwidth}
	\centering
	\includegraphics[width=\linewidth]{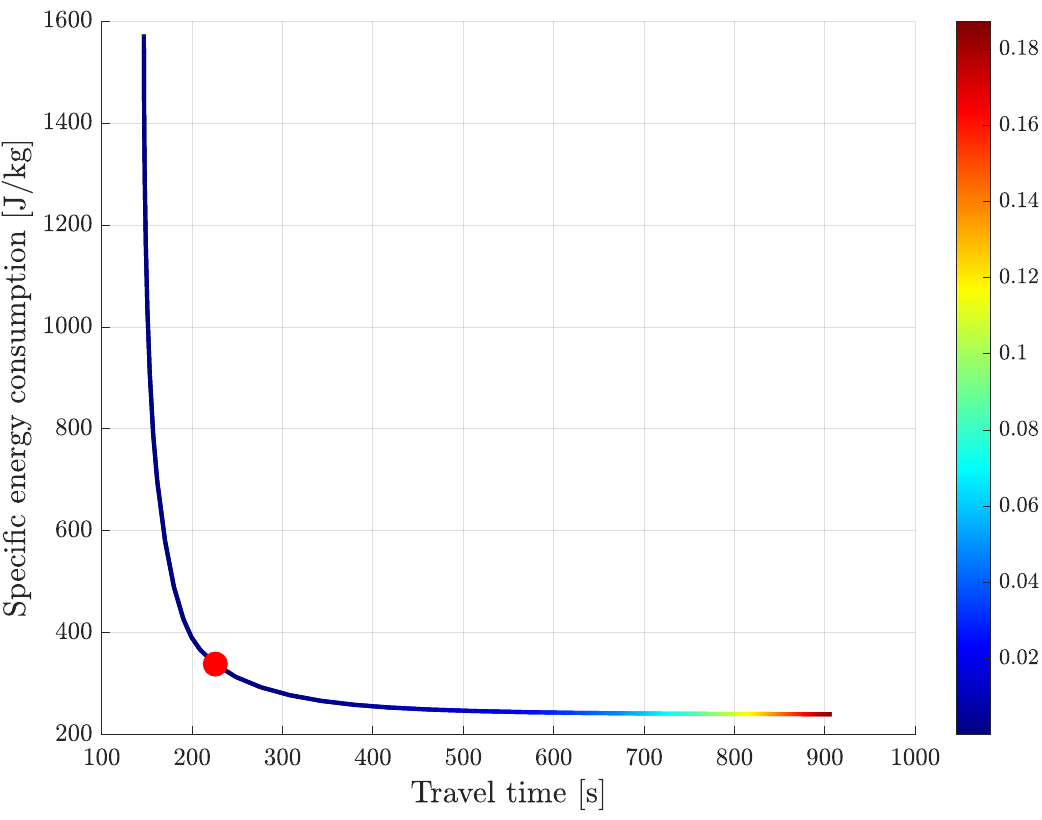}
	\caption{Monaco Grand Prix circuit.}
	\label{fig:MonacoPareto}
    \end{subfigure}
    \hfill
    \begin{subfigure}{0.49\columnwidth}
        \centering
        \includegraphics[width=\linewidth]{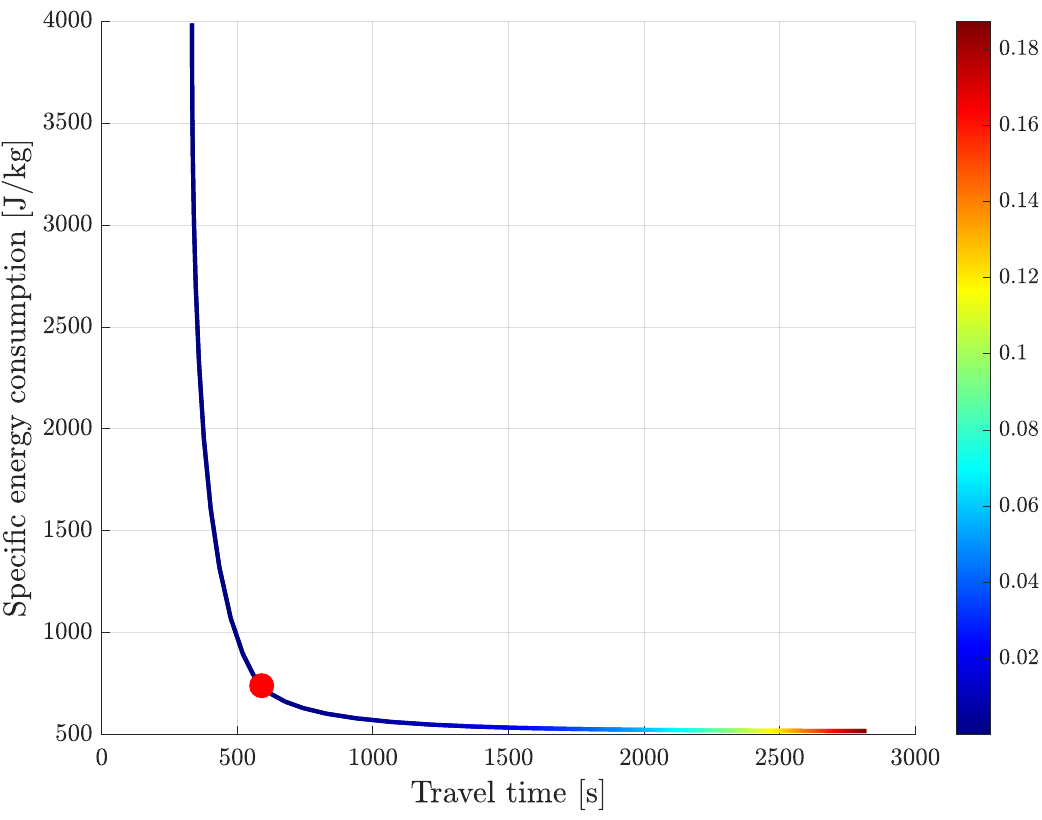}
        \caption{Urban path in Parma.}
        \label{fig:UniPRPareto}
    \end{subfigure}
	\begin{subfigure}{0.49\columnwidth}
	\centering
	\includegraphics[width=\linewidth]{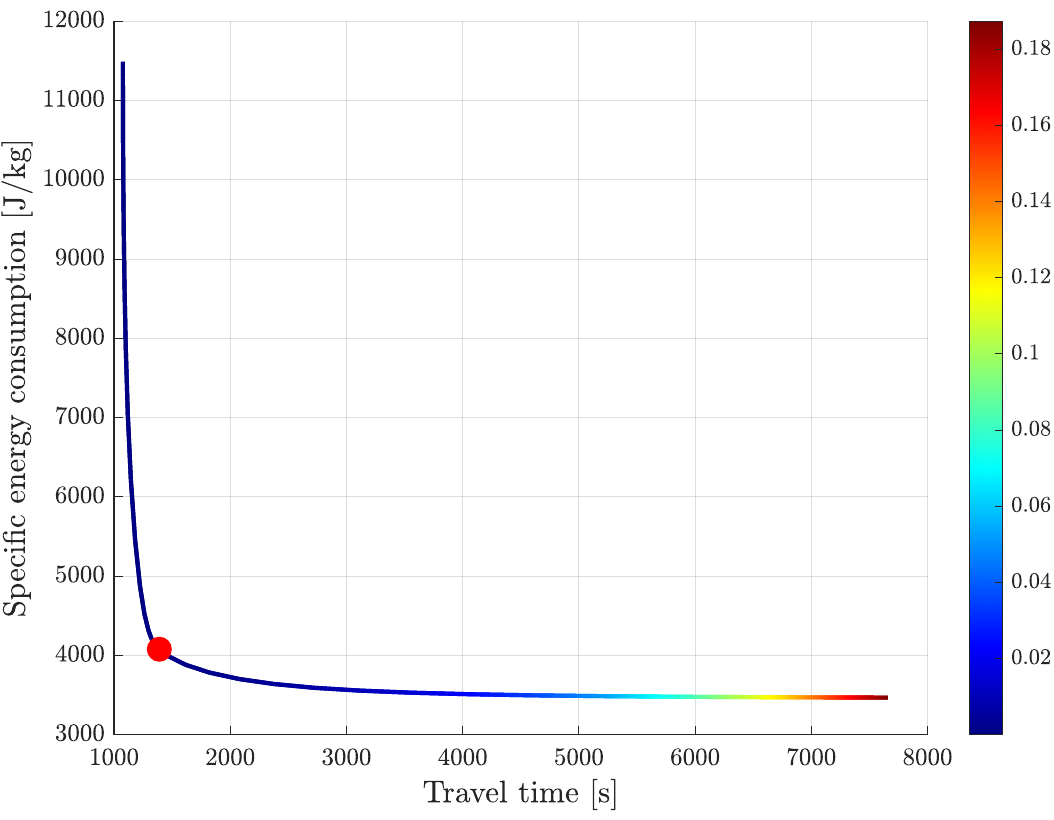}
	\caption{Passo Pordoi.}
	\label{fig:Arabba-CanazeiPareto}
    \end{subfigure}
    \hfill
    \begin{subfigure}{0.49\columnwidth}
        \centering
        \includegraphics[width=\linewidth]{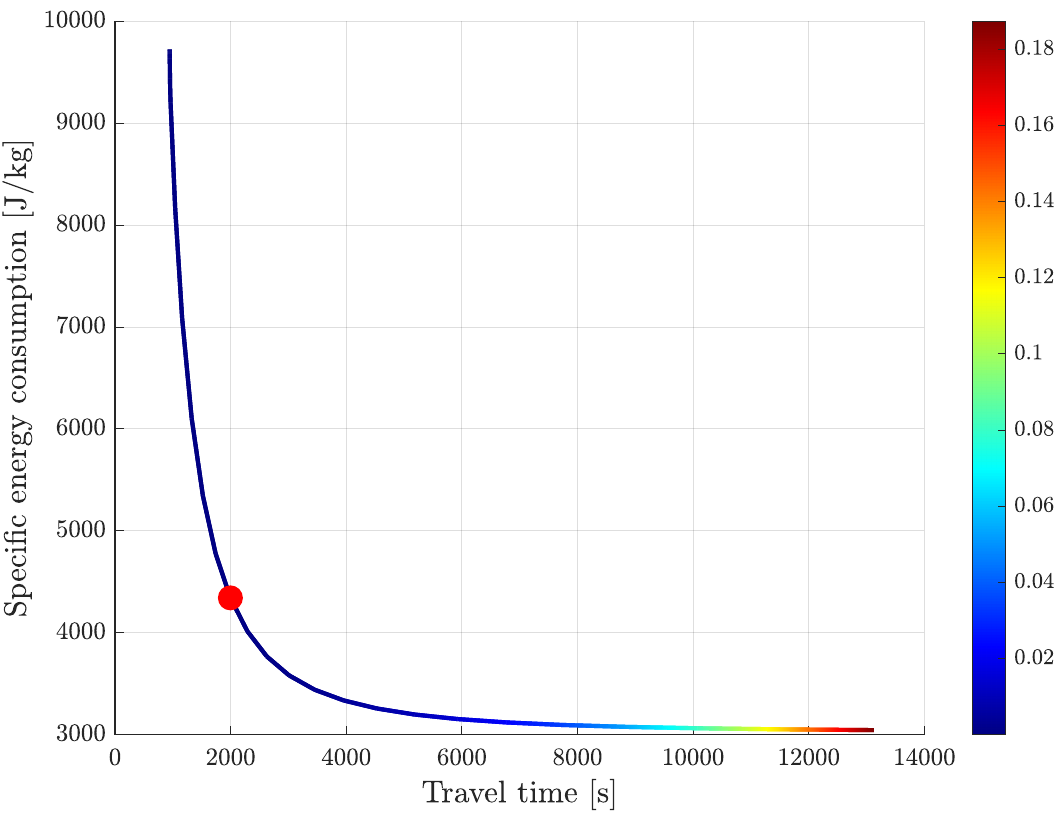}
        \caption{Section of Autostrada A1.}
        \label{fig:Rioveggio-BarberinoPareto}
    \end{subfigure}
    \caption{Pareto fronts for the four scenarios. The Pareto solutions corresponding to value $\lambda = 5\cdot10^{-4}$ are highlighted with a red circle.}
    \label{fig:Pareto_fronts}
\end{figure}
Recalling Remark~\ref{rem:pareto}, solving Problem~\eqref{eq:probfix} with different $\lambda$ values allows computing different Pareto optimal solutions
and, thus, allows representing the Pareto fronts for the four scenarios (see Figure~\ref{fig:Pareto_fronts}), which highlight the trade-offs between travel time and energy consumption. As shown in Figure~\ref{fig:Pareto_fronts}, the extreme values of $\lambda$, which correspond to the two opposite ends of the Pareto front, are undesirable: in one case, a slight improvement in travel time leads to a significant increase in energy consumption, while in the other, a marginal reduction in energy consumption results in a drastic increase in travel time.
Hence, a good practical compromise for typical road vehicles in everyday driving scenarios could be to consider values of $\lambda$ around $5\cdot10^{-4}$ (highlighted with a red circle in Figure~\ref{fig:Pareto_fronts}), which allow for good travel times with moderate energy consumption.
Indeed, much smaller values of $\lambda$ would slightly reduce travel times at the expenses of much larger energy consumptions and, viceversa, much larger values of $\lambda$ would lead towards minimum energy consumptions but with unacceptable high travel times.
In other words, choosing values of $\lambda$ away from $5\cdot10^{-4}$ cannot improve either of the two objective without detrimentally increasing the other.
Of course, the choice of $\lambda$ is context-dependent, if we were to plan speed for a racing car during a race, we would not reason in terms of good trade-offs between travel time and energy consumption but rather we would almost completely neglect the energetic term and choose values of $\lambda$ very close or equal to $0$ in order to privilege travel times.

\begin{figure}[h!]
	\centering
	\includegraphics[width=.6\columnwidth]{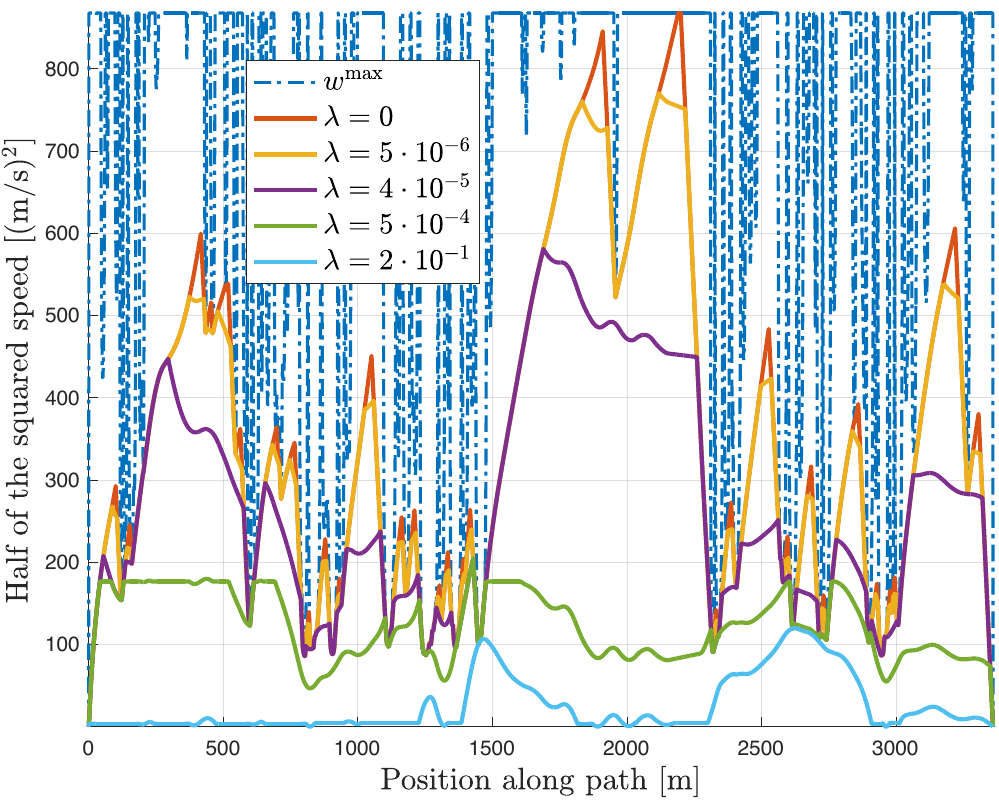}
	\caption{Some optimal half of the squared speed profiles for different values of $\lambda$.}
	\label{fig:MonacoProfiles}
\end{figure}

To better illustrate how parameter $\lambda$ influences the shape of speed profiles, in Figure~\ref{fig:MonacoProfiles}, we show some optimal half of the squared speed profiles corresponding to different values of $\lambda$ in the Monaco circuit. Similar profiles can be obtained for the other example scenarios. Figure~\ref{fig:MonacoProfiles} also shows the maximum half of the squared speed profile, which, due to the fact that the circuit is rather twisty, appears quite irregular.
As can be seen, for values of $\lambda$ close to 0, energy consumption is largely neglected. Consequently, the vehicle accelerates as aggressively as possible and breaks sharply to minimize travel time. However, as $\lambda$ increases, the ``spikes'' observed in the previous profile are progressively smoothed out, as reference speeds $\wrefp$ and $\wrefm$ decrease in response to the increasing weight assigned to energy consumption.
Finally, in Figure~\ref{fig:MonacoProfiles}, we also report the speed profile corresponding to $\lambda = 2\cdot10^{-1}$, which, although not practically meaningful, still exhibits a non-trivial shape due to the many uphill and downhill segments of the racing circuit.

\section{Conclusions}
\label{sec:concl}
In this paper we addressed the problem of planning the speed of a vehicle along an assigned path by minimizing a weighted sum of energy consumption and travel time.
The resulting mathematical model is a non-convex problem. We proved that the feasible region of this problem is a lattice which, if not empty, admits a minimum element $y$ and a maximum element $z$.
We introduced a feasibility-based bound-tightening technique which either allows to establish that the feasible region is empty, or it converges to the minimum and maximum elements of the feasible region.
Replacing the (non-convex) constraints of the non-convex problem with the lower limits $y$ and upper limits $z$ for the variables leads to a convex relaxation of the non-convex problem. We proved that such convex relaxation
is exact (i.e., that its optimal solution is also optimal for the original non-convex problem).
After studying the properties of optimal solutions of the convex relaxation, we exploited them to define a DP approach, which
is able to return a feasible solution of the problem. We proved that the optimal solution of the problem is close to at least one of the feasible solutions visited by the DP approach.
The complexity of such approach is quadratic with respect to the number of discretization steps $n$, thus allowing to return a solution very close to the optimal one with computing times much smaller than those required by well-known commercial solvers to solve the exact convex relaxation. Both the quadratic dependence on $n$ of the computing times, and the considerable computational gain with respect to commercial solvers
are confirmed through different computational experiments. 
Finally, we note that a MATLAB app built on the speed planning framework of this paper was presented in~\cite{IFAC-ACE}, as a tool for teaching multi-objective optimization problems in engineering education.

In the proposed model we employed specific definitions of the force and power constraints. We could consider 
the following, more general, definition of the minimum/maximum force and power constraints:
\begin{equation}
\label{eq:gencase}
\eta_i(w)\leq f_i\leq \beta_i(w_i) \qquad i \in \{1,\ldots,n-1\}.
\end{equation}
In the previous sections, we employed the following specific definitions:
\begin{gather*}
\beta_i(w)=\min\left\{g\mu\cos(\alpha_i),\frac{P_{\max}}{M\sqrt{2w}}\right\},\\
\eta_i(w)=-g\mu\cos(\alpha_i).
\end{gather*}
A possible future development is to prove that the results discussed in this paper can also be extended to more general functions $\beta_i$ and $\eta_i$, provided that these functions have a bounded variability and provided that the discretization step $h$ is small enough.

\bibliographystyle{plain} 
\bibliography{biblio}

\appendix
\section{Preliminary results for Theorem~\ref{thrm:main}}
\label{sec:prelim_res}
The first result deals with consecutive steps where the force is not null at an optimal solution
of Problem~\eqref{eq:probfix}.
\begin{prop}
\label{prop:simplest}
Consider an optimal solution $(w^\star,f^\star)$ of Problem~\eqref{eq:probfix} and assume that one of cases 
{\bf a)}--{\bf d)} of Proposition~\ref{prop:combposs} holds at some step $j$ (i.e., $f^\star_{j-1}f^\star_j\neq 0$). Then $w^\star_j\in W_j$. 
\end{prop}
\begin{proof}
The result immediately follows from Proposition~\ref{prop:combposs}.
\end{proof}
Next, we need the following proposition,
showing that we can bound from above the distance between the value of $w_j$ at an optimal solution of Problem~\eqref{eq:probfix} and one of the speeds in $W_j$ also when only one between $f_{j-1}$ and $f_j$ is different from $0$
(i.e., we are in cases {\bf e)}--{\bf h)} of Proposition~\ref{prop:combposs}).
\begin{prop}
\label{prop:distanceh}
Consider an optimal solution $(w^\star, f^\star)$ of Problem~\eqref{eq:probfix}, and assume that one of the cases {\bf e)}--{\bf h)} of Proposition~\ref{prop:combposs} holds at step $j$. 
Then, for some constant $C$ we have that:
\[
\min_{\bar{w}\in W_j} |\bar{w}-w^\star_j| \leq h C.
\]
\end{prop}
\begin{proof}
We only consider case {\bf f)} of Proposition~\ref{prop:combposs} with $f^\star_{j-1}>0, f^\star_j=0$. The proof for the other cases is analogous.
In this case we either have $w^\star_j=y_j$ or $w^\star_j\leq \wrefp$. In the former case, we obviously have $\min\limits_{\bar{w}\in W_j} |\bar{w}-w^\star_j|=0$, and we are done.

Otherwise, if $j=2$, then $w^\star_{j-1}=w^\star_1=y_1=z_1=\winit$, and we have that
\begin{align*}
|w^\star_2-z_2|&=|w^\star_2\!-\!w^\star_1\!+\!z_1\!-\!z_2|\leq |w^\star_2\!-\!w^\star_1|+|z_2\!-\!z_1|\leq\\
&\leq 2h\left[\gamma z_1 +g(|\sin(\alpha_1)| + (c+\mu)|\cos(\alpha_1)|)\right],
\end{align*}
which proves the result.

If $j>2$, then we also need to consider the sign of $f^\star_{j-2}$.
\begin{itemize}
\item If $f^\star_{j-2}<0$, in view of case {\bf b)} of Proposition~\ref{prop:combposs} with $f^\star_{j-2}<0, f^\star_{j-1}>0$, we must have $w^\star_{j-1}=z_{j-1}$.
Therefore, as before:
\begin{align*}
|w^\star_j-z_j| &= |w^\star_j-w^\star_{j-1}+z_{j-1}-z_j|
\leq |w^\star_j-w^\star_{j-1}|+|z_j-z_{j-1}|\leq  \\
&\leq 2h\left[\gamma z_{j-1} +g(|\sin(\alpha_{j-1})| + (c+\mu)|\cos(\alpha_{j-1})|)\right].
\end{align*}
\item If $f^\star_{j-2}>0$, in view of case {\bf d)} of Proposition~\ref{prop:combposs} with $f^\star_{j-2}>0, f^\star_{j-1}>0$, we have that $w^\star_{j-1}=z_{j-1}$ or $w^\star_{j-1}=y_{j-1}$ or $w^\star_{j-1}=\wrefp$.
The first two cases can be dealt with as before. If $w^\star_{j-1}=\wrefp$, we have that
\begin{align*}
|w^\star_j - w^\star_{j-1}| &=|w^\star_j-\wrefp|\leq\\
&\leq h\left[\gamma \wrefp +g(|\sin(\alpha_{j-1})| + (c+\mu)|\cos(\alpha_{j-1})|)\right].
\end{align*}
\item If $f^\star_{j-2}=0$, in view of case {\bf h)} of Proposition~\ref{prop:combposs} with $f^\star_{j-2}=0, f^\star_{j-1}>0$, we have that $w^\star_{j-1}=z_{j-1}$ or $w^\star_{j-1}\geq \wrefp$. The first case
can be dealt with as before.
In the second case we recall that $w^\star_j\leq \wrefp$, which together with $w^\star_{j-1}\geq \wrefp$, implies
\begin{align*}
|w^\star_j-\wrefp| & \leq |w^\star_j-w^\star_{j-1}| \\
& \leq h\left[\gamma w^\star_{j-1} +g(|\sin(\alpha_{j-1})| \!+\! (c + \mu)|\cos(\alpha_{j-1})|)\right]\!\leq \\
& \leq h\left[\gamma z_{j-1} +g(|\sin(\alpha_{j-1})| + (c+ \mu)|\cos(\alpha_{j-1})|)\right].
\end{align*}
\end{itemize}
\end{proof}
Finally, we need the following proposition, stating that the distance between the speeds when we perform a null-force move starting from
different initial speeds is bounded from above by the distance between such initial speeds.
\begin{prop}
\label{prop:distspeed}
Consider two null-force curves $w_{j}^{j,0}(w_j)$, \ldots, $w_{j+i}^{j,0}(w_j)$ and $w_{j}^{j,0}(\bar{w}_j),\ldots, w_{j+i}^{j,0}(\bar{w}_j)$ 
defined as in~\eqref{eq:nullforcurve}
with different initial speeds $w_j$ and $\bar{w}_j$.
Then:
\[
w_{j+r}^{j,0}(w_j)-w_{j+r}^{j,0}(\bar{w}_j)=(1-h\gamma)^r(w_j-\bar{w}_j),\ r \in \{0,\ldots,i\}.
\]
\end{prop}
\begin{proof}
The proof is by induction. The result is obviously true for $r=0$. Let us assume that it is true for $r-1$ and let us prove it holds also for $r$.
By definition~\eqref{eq:nullforcurve} we have that
\begin{align*}
w_{j+r}^{j,0}(w_j)=&(1-h\gamma)w_{j+r-1}^{j,0}(w_j) - h g(\sin(\alpha_{j+r-1}) + c\cos(\alpha_{j+r-1}))\\
w_{j+r}^{j,0}(\bar{w}_j)=&(1-h\gamma) w_{j+r-1}^{j,0}(\bar{w}_j) - h g(\sin(\alpha_{j+r-1}) + c\cos(\alpha_{j+r-1})).
\end{align*}
By taking the difference of these two equations we have that
\begin{align*}
w_{j+r}^{j,0}(w_j)-w_{j+r}^{j,0}(\bar{w}_j)
&=(1-h\gamma)\left[w_{j+r-1}^{j,0}(w_j)-w_{j+r-1}^{j,0}(\bar{w}_j)\right]=\\
&=(1-h\gamma)^r(w_j-\bar{w}_j),
\end{align*}
where the last equality follows from the inductive assumption and concludes the proof.
\end{proof}
Now, let us consider the subset $\bar{W}$ of the feasible region of the exact convex relaxation~\eqref{eq:newconvrel} made up of all $(w,f)\in \mathbb{R}^{2n-1}$ which fulfill the following:
\begin{gather*}
w_j\in W_j,\ w_{j+h}\in W_{j+h},\ w_{j+r}\not \in W_{j+r},\\
 r\in \{1,\ldots,h-1\},\ j\in \{1,\ldots,n-2\},\ h\geq 2\\ 
\Downarrow \\ 
f_{j}=\cdots=f_{j+h-2}=0.
\end{gather*}
Note that this is the set of speeds explored by the DP approach, since in that approach a speed $w_j$ different from one of the four speeds in $W_j$ can only be reached while moving along
a null-force move.
We remark that an optimal solution $(w^\star,f^\star)$ of~\eqref{eq:newconvrel} does not necessarily belong to $\bar{W}$ but we will prove below that it is close to it.

\section{Building a solution in $\bar{W}$ close to $(w^\star,f^\star)$}
\label{sec:solclose}  
We would like to show that, given an optimal solution 
$(w^\star,f^\star)$ of Problem~\eqref{eq:newconvrel}, we are able to build a solution $(\bar{w},\bar{f})\in \bar{W}$ close to it.
First of all, by Proposition~\ref{prop:combposs}
for each $j$ such that $f_{j-1}^\star,f_j^\star\neq 0$, we have that $w_j^\star\in W_j$, and we set $\bar{w}_j=w_j^\star$. Next, let us consider a sequence of forces with the following property for some $j, h\geq 2$:
\begin{equation}
\label{eq:seq}
f_{j-1}^\star,f_j^\star\neq 0,\ f_{j+1}^\star=\!\cdots\!=f_{j+h-1}^\star=0,\ f_{j+s}^\star,f_{j+s+1}^\star\neq 0.
\end{equation}
Note that we must have $\bar{w}_j=w_j^\star$ and $\bar{w}_{j+s+1}=w_{j+s+1}^\star$. If $f_j^\star>0$, we consider two possible options:
\[
\begin{array}{ll}
\mbox{{\tt Option A}}: & \bar{f}_j=0\implies \bar{w}_{j+1}<w_{j+1}^\star \\ 
\mbox{{\tt Option B}}: & \bar{f}_j=\min\left\{g\mu\cos(\alpha_{j}),\frac{\Pmax}{M\sqrt{2w_j^\star}}\right\} \implies \bar{w}_{j+1}\geq w_{j+1}^\star.
\end{array}
\]
If $f_j^\star<0$, the two options become:
\[
\begin{array}{ll}
\mbox{{\tt Option A}}: & \bar{f}_j=-g\mu\cos(\alpha_{j}) \implies \bar{w}_{j+1}\leq w_{j+1}^\star \\
\mbox{{\tt Option B}}: & \bar{f}_j=0\implies \bar{w}_{j+1}>w_{j+1}^\star.
\end{array}
\]
The choice between {\tt Option A} and {\tt Option B} depends on value $w_{j+s+1}^\star$. We proceed as follows:
\begin{equation}
\label{eq:choicopt}
\begin{aligned}
\mbox{Choose }&\mbox{{\tt Option A} if }  w_{j+s+1}^\star=y_{j+s+1}\\
& \mbox{ or }\ w_{j+s+1}^\star\in \{\wrefp, \wrefm\},\ w_{j+s+1}^\star<w_{j+s}^\star\\ 
\mbox{Choose }&\mbox{{\tt Option B} if }  w_{j+s+1}^\star=z_{j+s+1}\\
& \mbox{ or }\ w_{j+s+1}^\star\in \{\wrefp, \wrefm\},\ w_{j+s+1}^\star>w_{j+s}^\star.
\end{aligned}
\end{equation}
With these choices, also recalling Proposition~\ref{prop:distspeed}, a null-force move starting at $\bar{w}_{j+1}$ will meet the minimum speed profile curve $y$ (if $w_{j+s+1}^\star=y_{j+s+1}$) or the maximum speed profile curve $z$  (if  $w_{j+s+1}^\star=z_{j+s+1}$), or one of the two constant curves associated to reference speeds $\wrefp, \wrefm$  (if  $w_{j+s+1}^\star\in \{\wrefp, \wrefm\}$), at some step
$j+r$, with $r\leq s+1$.
We can refine~\eqref{eq:choicopt} by the following proposition.
\begin{prop}
\label{prop:opta}
If $w_{j+s+1}^\star=\wrefp$, then either $w_{j+s}^\star=z_{j+s}$, or $w_{j+s}^\star\geq w_{j+s+1}^\star$.
If $w_{j+s+1}^\star=\wrefm$, then either $w_{j+s}^\star=y_{j+s}$, or $w_{j+s}^\star\leq w_{j+s+1}^\star$.
\end{prop}
\begin{proof}
Assume that $w_{j+s+1}^\star=\wrefp$.
According to Proposition~\ref{prop:combposs} part {\bf d)}, we have that $w_{j+s+1}^\star=\wrefp$ implies $f_{j+s}^\star,f_{j+s+1}^\star>0$. Moreover, according to Proposition~\ref{prop:combposs} part {\bf h)},  
$f_{j+s-1}^\star=0,f_{j+s}^\star>0$ implies that $w_{j+s}^\star \geq  w_{j+s+1}^\star=\wrefp$ or that $w_{j+s}^\star=z_{j+s}$.

Similarly, if we assume that $w_{j+s+1}^\star=\wrefm$, then
according to Proposition~\ref{prop:combposs} part {\bf a)}, we have that $w_{j+s+1}^\star=\wrefm$ implies $f_{j+s}^\star,f_{j+s+1}^\star<0$. Moreover, according to Proposition~\ref{prop:combposs} part {\bf g)},  
$f_{j+s-1}^\star=0,f_{j+s}^\star>0$ implies that $w_{j+s}^\star \leq  w_{j+s+1}^\star=\wrefm$ or that $w_{j+s}^\star=y_{j+s}$.
\end{proof}
In view of this proposition, we can modify~\eqref{eq:choicopt} as follows:
\begin{equation}
\label{eq:choicopt1}
\begin{aligned}
\mbox{Choose {\tt Option A} if }&  w_{j+s+1}^\star=y_{j+s+1} \mbox{ or }\ w_{j+s+1}^\star=\wrefp\leq w_{j+s}^\star\\
\mbox{Choose {\tt Option B} if }&  w_{j+s+1}^\star=z_{j+s+1} \mbox{ or }\ w_{j+s+1}^\star=\wrefm\geq w_{j+s}^\star.
\end{aligned}
\end{equation}
Now, we define the values of  $\bar{w}$ and $\bar{f}$ depending on value $w_{j+s+1}^\star$. The following proposition deals with the case
$w_{j+s+1}^\star\in \{y_{j+s+1},z_{j+s+1}\}$.
Note that, according to notation~\eqref{eq:nullforcurve}, here and in what follows, $w_{j+r}^{j+r-1,0}(w_{j+r-1})$ denotes a one-step null-force move from step $j+r-1$ to step $j+r$, starting from half squared speed $w_{j+r-1}$, that is:
\[
w_{j+r}^{j+r-1,0}(w_{j+r-1})=(1-h\gamma)w_{j+r-1} - hg(\sin(\alpha_{j+r-1}) + c\cos(\alpha_{j+r-1})).
\]
\begin{prop}
\label{prop:caseyz}
Let $w_{j+s+1}^\star\in \{y_{j+s+1},z_{j+s+1}\}$. If we set, for $r \in \{2,\ldots,s\}$:
\begin{itemize}
\item
$\bar{w}_{j+r}=\min\{z_{j+r},w_{j+r}^{j+r-1,0}(\bar{w}_{j+r-1})\}$, if $w_{j+s+1}^\star=z_{j+s+1},$
\item
$\bar{w}_{j+r}=\max\{y_{j+r},w_{j+r}^{j+r-1,0}(\bar{w}_{j+r-1})\}$, if $w_{j+s+1}^\star=y_{j+s+1}$,
\end{itemize}
we have that $\bar{w}_{j+r}\geq w_{j+r}^\star$ for all $r\in \{1,\ldots,s+1\}$, if $w_{j+s+1}^\star=z_{j+s+1}$, and
$\bar{w}_{j+r}\leq w_{j+r}^\star$ for all $r\in\{1,\ldots,s+1\}$, if $w_{j+s+1}^\star=y_{j+s+1}$. Moreover, for $r \in \{1,\ldots,s+1\}$:
\[
|\bar{w}_{j+r}- w_{j+r}^\star|\leq |\bar{w}_{j+1}- w_{j+1}^\star|=O(h).
\]
\end{prop}
\begin{proof}
We only discuss the case $w_{j+s+1}^\star=z_{j+s+1}$. We prove the result by induction. Since in this case we are considering {\tt Option B}, we have that
$\bar{w}_{j+1}\geq w_{j+1}^\star$. Now, assume that $\bar{w}_{j+r-1}\geq w_{j+r-1}^\star$ and that $|\bar{w}_{j+r-1}- w_{j+r-1}^\star|\leq |\bar{w}_{j+1}- w_{j+1}^\star|$ for some $r\geq 2$.
We need to prove that the same holds at step $j+r$. If $\bar{w}_{j+r}=z_{j+r}$, then $z_{j+r}\geq w^\star_{j+r}$ implies that 
$\bar{w}_{j+r}\geq w^\star_{j+r}$. Instead, if $\bar{w}_{j+r}= w_{j+r}^{j+r-1,0}(\bar{w}_{j+r-1})$, in view of Proposition~\ref{prop:distspeed}, we have that
\[
 w_{j+r}^{j+r-1,0}(\bar{w}_{j+r-1})\geq w_{j+r}^{j+r-1,0}(w^\star_{j+r-1})=w_{j+r}^\star,
 \]
where the last equality follows from $f_{j+r-1}^\star=0$.
Moreover, again in view of Proposition~\ref{prop:distspeed} and by the inductive assumption, we have that
\begin{align*}
\bar{w}_{j+r}\!-\!w^\star_{j+r}&\leq w_{j+r}^{j+r-1,0}(\bar{w}_{j+r-1})\!-\!w_{j+r}^{j+r-1,0}(w^\star_{j+r-1})\leq\\
&\leq|\bar{w}_{j+r-1}- w_{j+r-1}^\star|\leq  |\bar{w}_{j+1}- w_{j+1}^\star|.
\end{align*}
\end{proof}
Notice that at all steps where $\bar{w}_{j+r}=w_{j+r}^{j+r-1,0}(\bar{w}_{j+r-1})$, we have that $\bar{f}_{j+r-1}=f^\star_{j+r-1}=0$.

Next, we discuss the cases in which $w_{j+s+1}^\star=\wrefp$ (the cases $w_{j+s+1}^\star=\wrefm$ can be dealt with in a completely analogous way).
The definition of $\bar{w}$ is more complicated with respect to the previous case.
As already observed, we must have that $w_{j+s}^\star \geq  w_{j+s+1}^\star=\wrefp$ or, alternatively, that $w_{j+s}^\star=z_{j+s}$. The last case can be addressed as already done in
Proposition~\ref{prop:caseyz}. Therefore, we only consider the case $w_{j+s}^\star \geq  w_{j+s+1}^\star$.
In view of Proposition~\ref{prop:distspeed} and the fact that we need to consider {\tt Option A}, that is, $w_{j+1}^\star\geq \bar{w}_{j+1}$, we must have that
\[
w_{j+s+1}^{j+1,0}(\bar{w}_{j+1})\leq w_{j+s+1}^{j+1,0}(w^\star_{j+1}).
\]
Now, let us introduce the following interval:
\begin{align*}
J(\bar{w}_{j+r-1}) =&\left[\min\left\{\bar{w}_{j+r-1},w_{j+r}^{j+r-1,0}(\bar{w}_{j+r-1})\right\}, \right.\\
&\left.\max\left\{\bar{w}_{j+r-1},w_{j+r}^{j+r-1,0}(\bar{w}_{j+r-1})\right\}\right],
\end{align*}
and, for $r \in \{2,\ldots,s\}$, let us introduce the following sets:
\[
\Omega(\bar{w}_{j+r-1})=
\begin{cases}
\{\wrefp,\ w_{j+r}^{j+r-1,0}(\bar{w}_{j+r-1})\}, & \mbox{if } \wrefp\in J(\bar{w}_{j+r-1}) \\
\{ w_{j+r}^{j+r-1,0}(\bar{w}_{j+r-1})\}, & \mbox{otherwise.}
\end{cases}
\]
In other words, we are distinguishing between the steps where a one-step null force move crosses constant line $\wrefp$, that is, those for which $\wrefp\in J(\bar{w}_{j+r-1}) $, and all the other steps.
Let us define values $\bar{w}$ as follows:
\begin{equation}
\label{eq:barwwref}
\bar{w}_{j+r}=
\begin{cases}
\wrefp, &\mbox{if } (w_{j+r-1}^\star-\wrefp)(w_{j+r}^\star-\wrefp)<0 \\
\underset{{w\in \Omega(\bar{w}_{j+r-1})}}{\arg\min} |w-w^\star_{j+r}|, &\mbox{otherwise}.
\end{cases}
\end{equation}
We can prove the following proposition, analogous to Proposition~\ref{prop:caseyz}, for the case when $w^\star_{j+s+1}=\wrefp$.
\begin{prop}
Let $w_{j+s+1}^\star=\wrefp\leq w_{j+s}^\star$. 
Let $\bar{w}$ be defined as in~\eqref{eq:barwwref}. Then, 
\[
|\bar{w}_{j+r}- w_{j+r}^\star|=O(h) \qquad r \in \{1,\ldots,s\}.
\]
\end{prop}
\begin{proof}
The result is true for $r=1$ in view of the definition of $\bar{w}_{j+1}$ following {\tt Option A}. Therefore, we need to prove that in case it is true for some $r-1$, it also holds for $r$, for any $r\geq 2$.
If $\bar{w}_{j+r}=w_{j+r}^{j+r-1,0}(\bar{w}_{j+r-1})$, recalling that $w_{j+r}^\star=w_{j+r}^{j+r-1,0}(w^\star_{j+r-1})$, in view of Proposition~\ref{prop:distspeed}, we have that:
\begin{align*}
\bar{w}_{j+r}-w^\star_{j+r}&=w_{j+r}^{j+r-1,0}(\bar{w}_{j+r-1})-w_{j+r}^{j+r-1,0}(w^\star_{j+r-1})=\\
&=(1-h\gamma)(\bar{w}_{j+r-1}-w_{j+r-1}^\star),
\end{align*}
and the results follows from the inductive assumption.
Therefore, we only need to prove the result when $\bar{w}_{j+r}=\wrefp$. We first consider $r$ such that $ (w_{j+r-1}^\star-\wrefp)(w_{j+r}^\star-\wrefp)\geq 0$.
But in such case, by definition of $\bar{w}_{j+r}$ we have that
\begin{align*}
|\bar{w}_{j+r}\!-\!w_{j+r}^\star|&\!\leq\! |w_{j+r}^{j+r-1,0}(\bar{w}_{j+r-1})\!-\!w_{j+r}^{j+r-1,0}(w^\star_{j+r-1})|=\\
&\!=\!|(1-h\gamma)(\bar{w}_{j+r-1}\!-\!w_{j+r-1}^\star)|,
\end{align*}
form which the result follows again by the inductive assumption.
Finally, we consider the case $\bar{w}_{j+r}=\wrefp$ and $r$ such that $ (w_{j+r-1}^\star-\wrefp)(w_{j+r}^\star-\wrefp)< 0$. Assume that
$w_{j+r-1}^\star>\wrefp, w_{j+r}^\star<\wrefp$. Then, $w_{j+r}^\star=w_{j+r}^{j+r-1,0}(w^\star_{j+r-1})> w_{j+r}^{j+r-1,0}(\wrefp)$, and:
\[
|\bar{w}_{j+r}-w_{j+r}^\star|= \wrefp-w_{j+r}^{j+r-1,0}(w^\star_{j+r-1})\leq
\wrefp-w_{j+r}^{j+r-1,0}(\wrefp)=O(h).
\]
Similarly, if 
$w_{j+r-1}^\star<\wrefp, w_{j+r}^\star>\wrefp$, Then $w_{j+r}^\star=w_{j+r}^{j+r-1,0}(w^\star_{j+r-1})< w_{j+r}^{j+r-1,0}(\wrefp)$, and:
\[
|\bar{w}_{j+r}-w_{j+r}^\star| = w_{j+r}^{j+r-1,0}(w^\star_{j+r-1})-\wrefp\leq
w_{j+r}^{j+r-1,0}(\wrefp)-\wrefp=O(h).
\]
\end{proof}
Note that in all cases in~\eqref{eq:barwwref} where $\bar{w}_{j+r} = w_{j+r}^{j+r-1,0}(\bar{w}_{j+r-1})$, we also have that $\bar{f}_{j+r-1}=f^\star_{j+r-1}=0$.

We still need to address some special cases. Indeed, sequences~\eqref{eq:seq} do not cover all possible cases since it may happen that for some index $s$:
\begin{equation}
\label{eq:subseq}
f_{j+s-1}^\star=0,\qquad f_{j+s}^\star\neq 0,\qquad f_{j+s+1}^\star=0.
\end{equation}
However, all sequences containing subsequences~\eqref{eq:subseq} can be converted into sequences of form~\eqref{eq:seq}. 
We only discuss the case $f_{j+s}^\star> 0$. In view of parts {\bf h)} and {\bf f)} of Proposition~\ref{prop:combposs}, we have that:
\begin{gather*}
w^\star_{j+s}\geq \wrefp\ \mbox{\tt or }\ w^\star_{j+s}=z_{j+s} \\
w^\star_{j+s+1}=y_{j+s+1}\  \mbox{\tt or }\ w^\star_{j+s+1}\leq \wrefp.
\end{gather*}
If $w^\star_{j+s}=z_{j+s}$, then we deal with subsequence~\eqref{eq:subseq} as if it were
\[
f_{j+s-1}^\star<0,\qquad f_{j+s}^\star> 0,\qquad f_{j+s+1}^\star=0,
\]
so that we are back to sequences~\eqref{eq:seq}.
If  $w^\star_{j+s+1}=y_{j+s+1}$, then we deal with subsequence~\eqref{eq:subseq} as if it were
\[
f_{j+s-1}^\star=0,\qquad f_{j+s}^\star> 0,\qquad f_{j+s+1}^\star<0,
\]
and again we are back to sequences~\eqref{eq:seq}.
Finally, if $w^\star_{j+s}\geq \wrefp \geq w^\star_{j+s+1}$, we fix $\bar{w}_{j+s}=\wrefp$ and we deal with subsequence~\eqref{eq:subseq} as if it were
\[
f_{j+s-1}^\star>0,\qquad f_{j+s}^\star> 0,\qquad f_{j+s+1}^\star=0.
\]

\section{Proof of Theorem~\ref{thrm:main}}
\label{sec:final_proof}
\begin{proof}
The fact that $|\bar{w}_j-w_j^\star|=O(h)$ for all $j\in\{1,\ldots,n\}$ follows from the discussion in Appendix~\ref{sec:solclose}.
Next, note that:
\[
| hF(\bar{w},\bar{f})-hF(w^\star,f^\star)| \leq \left|h F_1(w^\star)-hF_1(\bar{w})\right|+\left| hF_2(\bar{f})-hF_2(f^\star)\right|.
\]
If $w_j^\star>0$ for all $j\in {\mathcal{J}}_{\bar w}$, then the derivative of $F_1$ at each such point $w_j^\star$ is bounded and, consequently, we have that:
\[
\left|h F_1(w^\star)-hF_1(\bar{w})\right|=O(h).
\]
Finally, we observe that 
\[
hF_2(\bar{f})-hF_2(f^\star)=
h \eta\gamma \lambda M \sum_{j\in {\mathcal{J}}_f} \left(\max\{\bar{f}_j,0\}-\max\{f_j^\star,0\}\right),
\]
and the result follows after observing that:
\[
\left| \max\{f_j^\star,0\}-\max\{\bar{f}_j,0\} \right|\leq g \mu |\cos(\alpha_j)|\leq g\mu.
\]
\end{proof}

\end{document}